\documentclass[openany, amssymb, psamsfonts]{amsart}

\usepackage[margin=1.5in]{geometry}

\usepackage{macros}
\usepackage{mathbbol}
\usepackage{mathrsfs,comment}
\usepackage[usenames,dvipsnames]{color}
\usepackage[normalem]{ulem}
\usepackage{url}
\usepackage[all,arc,2cell]{xy}
\UseAllTwocells
\usepackage{enumerate}
\usepackage{enumitem}
\usepackage{hyperref}
\hypersetup{%
  bookmarksnumbered=true,%
  bookmarks=true,%
  colorlinks=true,%
  linkcolor=blue,%
  citecolor=blue,%
  filecolor=blue,%
  menucolor=blue,%
  pagecolor=blue,%
  urlcolor=blue,%
  pdfnewwindow=true,%
  pdfstartview=FitBH}

\usepackage{cleveref}
\usepackage{todonotes}
\usepackage{graphicx}
\usepackage{adjustbox}
\usepackage{amssymb}
\usepackage{tikz-cd}
\usepackage[backend=bibtex]{biblatex}
\usepackage{amsthm}
\def\makeautorefname#1#2{\expandafter\def\csname#1autorefname\endcsname{#2}}
\def\equationautorefname~#1\null{(#1)\null}
\makeautorefname{footnote}{footnote}%
\makeautorefname{item}{item}%
\makeautorefname{figure}{Figure}%
\makeautorefname{table}{Table}%
\makeautorefname{part}{Part}%
\makeautorefname{appendix}{Appendix}%
\makeautorefname{chapter}{Chapter}%
\makeautorefname{section}{Section}%
\makeautorefname{subsection}{Section}%
\makeautorefname{subsubsection}{Section}%
\makeautorefname{theorem}{Theorem}%
\makeautorefname{thm}{Theorem}%
\makeautorefname{cor}{Corollary}%
\makeautorefname{lem}{Lemma}%
\makeautorefname{prop}{Proposition}%
\makeautorefname{pro}{Property}
\makeautorefname{conj}{Conjecture}%
\makeautorefname{defn}{Definition}%
\makeautorefname{notn}{Notation}
\makeautorefname{notns}{Notations}
\makeautorefname{rem}{Remark}%
\makeautorefname{quest}{Question}%
\makeautorefname{exmp}{Example}%
\makeautorefname{ax}{Axiom}%
\makeautorefname{claim}{Claim}%
\makeautorefname{ass}{Assumption}%
\makeautorefname{asss}{Assumptions}%
\makeautorefname{con}{Construction}%
\makeautorefname{prob}{Problem}%
\makeautorefname{warn}{Warning}%
\makeautorefname{obs}{Observation}%
\makeautorefname{conv}{Convention}%

\theoremstyle{definition}
\newtheorem{thm}{Theorem}[section]
\newtheorem{cor}{Corollary}[section]
\newtheorem{prop}{Proposition}[section]
\newtheorem{lem}{Lemma}[section]

\newtheorem{defn}{Definition}[section]

\newtheorem{exmp}{Example}[section]
\newtheorem{notn}{Notation}[section]

\newtheorem{rem}{Remark}[section]

\makeatletter
\let\c@obs=\c@thm
\let\c@cor=\c@thm
\let\c@prop=\c@thm
\let\c@lem=\c@thm
\let\c@prob=\c@thm
\let\c@con=\c@thm
\let\c@conj=\c@thm
\let\c@defn=\c@thm
\let\c@notn=\c@thm
\let\c@notns=\c@thm
\let\c@exmp=\c@thm
\let\c@ax=\c@thm
\let\c@pro=\c@thm
\let\c@ass=\c@thm
\let\c@warn=\c@thm
\let\c@rem=\c@thm
\let\c@sch=\c@thm
\let\c@equation\c@thm
\numberwithin{equation}{section}
\makeatother

\title{Graded Monoidal Categories and Internalization}

\author{Bryce Goldman}

\begin{document}

\begin{abstract}
    In stable homotopy theory, there are a variety of ways to construct a point-set level smash product of spectra, each with its own advantages and limitations. These constructions often involve internalizing some graded external product. In this paper, we aim to formalize this process through a general categorical framework, providing a common language for these constructions. We propose a general theory for internalizing a graded monoidal product, investigate its 2-categorical properties, and use this framework to recover both Day convolution and the smash product of $S$-modules.
\end{abstract}

\maketitle

\tableofcontents

\section{Introduction}
\label{section:intro}

In homotopy theory one is often concerned only with objects of interest up to homotopy. However, when possible, one may wish to lift a homotopical phenomenon to the point-set level to be studied prior to passing to the homotopy category. For instance, one might desire a well-behaved smash product of spectra that lifts the symmetric monoidal structure of the stable homotopy category. However, in \cite{lewis}, it is shown that no such lift can exist without surrendering at least one property desirable of this product.

In fact, many weaker smash products of spectra exist, varying in which of these properties are compromised and what point-set model of spectra is favored. For example, one can construct a smash product of symmetric or orthogonal spectra using Day convolution \cite{MMSS, schwede, HSS}. Alternatively, if one wishes to use the coordinate-free spectra of \cite{EKMM} then restricting to $S$-modules yields another smash product. Both of these constructions entail some external graded product being pulled into a single degree to yield the desired internal smash product. In this paper, we propose a general categorical axiomatization of this procedure of internalizing a graded product: we define an additional structure on \emph{graded monoidal categories} which permits an internalization, show that the resulting product is monoidal, and recover the constructions of several classical smash products of spectra using this framework. 

We make the following contributions.
\begin{itemize}[label=\textendash, noitemsep, topsep=2pt]
    \item The structure of an internalization on a graded monoidal category, forming what we call an \emph{internalized category} (\Cref{def:internalization}).
    \item A framework for constructing an ordinary monoidal product on the degree one component of an internalized category, which can be equipped with additional structure allowing for the recovery of standard properties like symmetry and closure (\Cref{thm:internalization}, \Cref{cor:symmetric_internalization}, \Cref{cor:closed_internalization}).
    \item A 2-categorical theory of graded monoidal categories and internalization, including the corresponding notions of functors and natural transformations between them (\Cref{section:categorical_properties}).
    \item An outline of free-forgetful adjunction between graded monoidal categories and internalized categories, as well as the coreflective relationship between ordinary categories and graded monoidal categories, and between monoidal categories and internalized categories (\Cref{cor:gr-coreflective-cat}, \Cref{thm:int-section}, \Cref{thm:free-forget-adj-mon-int}).
    \item Derivations of Day convolution and the smash product of $S$-modules within the framework of internalization, providing a common language across both constructions using their graded monoidal structure (\Cref{prop:day_is_internalization}, \Cref{thm:wedgeS_internal}).
\end{itemize}

\subsection{Notation}
\label{subsection:notation}

Before we begin assembling the machinery of internalization, we make a few brief notes on the notation used throughout this paper. 

We denote by $\mathbb{1}$ the terminal category, consisting of a unique object and its identity morphism. Given a category $\C$ and objects $x, y \in \C$, we write $\C(x, y)$ for the set\footnote{For simplicity, we assume foundations that allow for all categories considered to be made locally small.} of morphisms with domain $x$ and codomain $y$. If $\C$ is instead enriched over some monoidal category $\V$, then in most cases we will write $\underline{\C}(x, y)$ for the hom-object in $\V$ of morphisms from $x$ to $y$, to emphasize the enrichment.

Additionally, we will write $\NN$ for the set of natural numbers (starting at $0$). Given $m \in \NN$ we denote by $\C^m$ the $m\Th$ power of $\C$, that is, the iterated product
\begin{equation*}
    \C^m := \prod_{i = 1}^m \C.
\end{equation*} Note that $\C^0 = \mathbb{1}$.

Finally, much of the content of the following sections involves objects which are $\NN$-graded. Consequently, a great deal of indices are required to keep track of degrees. When possible, we omit these indices, provided that the degree of the objects we are concerned with are either clear from context or unimportant. In particular, our diagrams will almost always be free of indices.

\section{External and Internal Products}
\label{section:external_and_internal}

Our aim is to formalize the procedure of converting an external, graded product into an ordinary monoidal product. To achieve this, we must first make sense of what it means for a (graded) category to be equipped with such an external product. We fix the notion of a monoidal product graded over the naturals in \Cref{subsection:graded_monoidal}, and introduce generalizations of a few common properties enjoyed by ordinary monoidal categories (namely, symmetry and closure). We then go on to make rigorous the structure which enables internalizing a graded product, and establish some basic results regarding internalization, in \Cref{subsection:internalization}.

\begin{rem}\label{rem:enriched} 
    All of the definitions and constructions in \Cref{subsection:graded_monoidal} and \Cref{subsection:internalization} can be interpreted more generally in the enriched setting. In particular, we could fix a finitely-complete\footnote{The completeness hypothesis is necessary to form the binary product of $\V$-categories.} closed monoidal category $(\V, \otimes_\V, I_\V)$ and consider external products (as $\V$-functors) on graded $\V$-categories. We omit a technical description of this generalization, as the interpretations of the arguments above would be at best cumbersome, and at worst entirely unreadable. 
    
    In \Cref{subsection:day_convolution} the internalized categories we are interested in will be enriched, so we note that this is possible only so that the reader is not completely blindsided by the sudden use of enriched categories.
\end{rem}

\subsection{Graded Monoidal Categories}
\label{subsection:graded_monoidal}

We begin by making the following definition:

\begin{defn} \label{def:graded_monoidal}
A \emph{graded monoidal category} $(\C_\smallbullet, \mu, 1)$ consists of the following data.
\begin{enumerate}
    \item For each $m \in \NN$ a category $\C_m$, called the \emph{degree $m$ component} of $\C_\smallbullet$. (We will also occasionally write $\C$ for the coproduct of the components of $\C_\smallbullet$.)
    \item A bifunctor $\mu_{m, n} \colon \C_m \times \C_n \to \C_{m + n}$ for each pair $m, n \in \NN$. We will often abuse notation by writing $\mu$ for the gestalt of these functors. (We will often omit degree indices; since $\mu$ is degree-additive, the degree of $\mu(x, y)$ is uniquely determined by the degrees of $x$ and $y$, so this is unambiguous.)
    \item An object $1 \in \C_0$, called the \emph{unit} of $\mu$.
    \item A natural isomorphism $\alpha_{m, n, r}$ for each $m, n, r \in \NN$, called the \emph{associators} of $\mu$:
    \begin{equation*}
        \begin{tikzcd}[column sep = 2.5em, row sep = 2.5em]
            \C_m \times \C_n \times \C_r \ar[r, "\C_m \times \mu_{n, r}"] \ar[d, "\mu_{m, n} \times \C_r"'] & \C_m \times \C_{n + r} \ar[d, "\mu_{m, n + r}"] \ar[dl, Rightarrow, shorten >= 4ex, shorten <= 3.5ex, "\alpha_{m, n, r}"']\\
            \C_{m + n} \times \C_r \ar[r, "\mu_{m + n, r}"'] & \C_{m + n + r} \period
        \end{tikzcd}
    \end{equation*}
    \item Natural isomorphisms $\lambda_m $ and $\rho_m$ for each $m \in \NN$, called the \emph{left unitors} and \emph{right unitors} of $\mu$, respectively:
    \begin{equation*}
        \begin{tikzcd}[column sep = 1.5em]
            & \C_0 \times \C_m \ar[dr, "\mu_{0,m}"], \ar[d, Rightarrow, shorten <= 1.5ex, shorten >= 1.5ex, "\lambda_m"] \\
            \C_m \ar[ur, "{(1, \C_m)}"] \ar[rr, equals] & {} & \C_m
        \end{tikzcd} \hspace{1.5em} 
        \begin{tikzcd}[column sep = 1.5em]
            C_m \ar[rr, equals] \ar[dr, "{(\C_m, 1)}"'] & {} & \C_m\\
            & \C_m \times \C_0 \ar[ur, "\mu_{m, 0}"'] \ar[u, Rightarrow, shorten <= 1.5ex, shorten >= 0.5ex, "\rho_m"'] \period
        \end{tikzcd}
    \end{equation*}
\end{enumerate} This data is further subject to the same coherence conditions as ordinary monoidal categories; that is, for each $x, y, z, w \in \C$ the following diagrams commute:
\begin{equation} \label{dgm:unity}
    \begin{tikzcd}[column sep = 2em, row sep = 2.5em]
        & \mu(x, y) &\\
        \mu(x, \mu(1, y)) \ar[ur, "{\mu(x, \lambda_y)}"] \ar[rr, "{\alpha_{x, 1, y}}"'] && \mu(\mu(x, 1), y) \ar[ul, "{\mu(\rho_x, y)}"']
    \end{tikzcd}
\end{equation}

\begin{equation} \label{dgm:associativity}
    \begin{tikzcd}[column sep = 0em, row sep = 2.5em]
        \mu(x, \mu(y, \mu(z, w))) \ar[rr, "{\alpha_{x, y, \mu(z, w)}}"] \ar[d, "{\mu(x, \alpha_{y, z, w})}"'] && \mu(\mu(x, y), \mu(z, w)) \ar[d, "{\alpha_{\mu(x, y), z, w}}"]\\
        \mu(x, \mu(\mu(y, z), w)) \ar[dr, "{\alpha_{x, \mu(y, z), w}}"'] && \mu(\mu(\mu(x, y), z), w)\\
        & \mu(\mu(x, \mu(y, z)), w) \ar[ur, "{\mu(\alpha_{x, y, z}, w)}"'] \period
    \end{tikzcd}
\end{equation}
\end{defn}

We think of the collection of categories $\C_\smallbullet$ as a single monoidal category, with the additional structure of a grading. Indeed, $\mu$ defines an ordinary monoidal structure on $\C := \bigsqcup_{m \in \NN} \C_m$, since each morphism of $\C$ is contained in some graded component $\C_m$ (so the naturality and coherence of of the structure maps can be checked degree-wise).

\begin{rem}[Coherence]
    \Cref{def:graded_monoidal} is a direct extension of the definition of an ordinary monoidal category: all of the coherence conditions are exactly the same as usual, ignoring the minor technicality of tracking degrees. Mac Lane's Coherence Theorem for monoidal categories \cite[Section VII]{maclane} still holds for a graded monoidal category $\C_\sbull$, since it holds for the monoidal category $\bigsqcup_{m \in \NN} \C_m$.
\end{rem}

\begin{rem}
It makes sense to allow more general notions of grading. In our setting, \Cref{def:graded_monoidal} makes sense for any collection of categories graded by a monoid $(M, \cdot, e),$ where the monoidal product has signature $\mu \colon \C_x \times \C_y \to \C_{x\cdot y}$ and the identity $1^\C$ lives in $\C_e$. One can further generalize to grading by a monoidal \emph{category}. For reasons that will become apparent in \Cref{section:external_and_internal} we are specifically interested in the case when the structure indexing a graded monoidal category is generated by a distinguished element (e.g. $1 \in \mathbb{N}$); for simplicity, we restrict ourselves to grading all of our categories by $\mathbb{N}$.
\end{rem}

\begin{rem}
Prior work has made use of related notions of monoidal categories equipped with some kind of grading. An earlier notion of graded monoidal category, aimed at abstracting certain equivariant methods, are graded by the delooping of a group \cite{frohlich-wall}, while more recent work defines a notion of grading by partial commutative monoids which that generalizes our specialization to $\NN$ \cite{earnshaw-nester-roman}. There is also an established body of work surrounding \emph{locally} graded categories over a monoidal category $\mathcal{V}$, which unify $\mathcal{V}$-enriched categories and $V$-actegories \cite{lucyshyn-wright, wood}. Our graded monoidal categories can also be viewed as instances of these by taking $\mathcal{V}$ to be $(\NN, +, 0)$ and requiring that no object is the domain or codomain of morphisms with different gradings.
\end{rem}

Examples of graded monoidal categories are abundant in nature:

\begin{exmp} 
    Any monoidal category $\C$ can be viewed as a graded monoidal category in a number of ways. One can form the \emph{trivial graded category} $\C(0)_\sbull$, which is given by $\C$ in degree 0, and which is the empty category in all higher degrees. (One could equally take $\C(0)_m$ to be the terminal category $\mathbb{1}$ with one object and one morphism for each $m > 0$.) The monoidal product on $\C$ then naturally extends to a graded monoidal product on $\C(0)_\sbull$. 

    The ordinary monoidal category $\C$ also gives rise to the \emph{constant graded category} $\overline{C}_\sbull$, defined by $\overline{\C}_m := \C$ for each $m \in \NN$. Again, the monoidal product on $\C$ extends to a graded monoidal product on $\overline{\C}_\sbull$.

    Finally, one can produce a \emph{free graded monoidal category} $\Free(\C)_\sbull$ on any given category $\C$, monoidal or not. We let $\Free(\C)_m := \C^m$ for each $m \in \NN$, where $\C^0 = \mathbb{1}$, and let $\Free(\C)_m \times \Free(\C)_n \to \Free(\C)_{m + n}$ be given by concatenation. This construction enjoys a handful of nice properties when $\C$ is monoidal, which are outlined in \Cref{section:categorical_properties}.
\end{exmp}

\begin{exmp}
    Every graded ring $S = \bigoplus_{d \geq 0} S_d$ induces a graded monoidal category $\mathcal{S}_\sbull,$ whose degree $m$ component is the cartesian product $\prod_{d \leq m} \categ{B} S_d$, where $\categ{B} S_d$ is the one-object category whose morphisms are the elements of the abelian groups $S_d$. The graded monoidal product is supplied by the multiplication of $S$. For instance, if $S = R[x]$ is a polynomial ring over some coefficient ring $R$ then $\mathcal{S}_m$ has a single object, with polynomials of degree $m$ as morphisms and addition as composition.
\end{exmp}

\begin{exmp} 
    Recall that if $M$ and $N$ are connected smooth manifolds of dimensions $m$ and $n$ (respectively), then the cartesian product $M \times N$ is a connected smooth $(m + n)$-manifold. The category $\categ{FinMan}$ of finite-dimensional connected smooth manifolds can be turned into a graded category $\categ{FinMan}_\sbull$, where $\categ{FinMan}_m$ is the full subcategory of $m$-manifolds. The categorical product of $\categ{FinMan}$ then forms a graded monoidal product on $\categ{FinMan}_\sbull$, since it is dimension-additive.
    
    Note that we have lost the data of morphisms between manifolds of distinct dimension in forming $\categ{FinMan}_\sbull$. As we will see in the following sections, the purpose of a graded monoidal category will often be to study the degree 1 component in particular. In practice, we will typically build the higher degrees of a graded monoidal category from degree 1, as opposed to stripping a category of enough data to reduce it to a graded monoidal category.
\end{exmp}

As a kind of specialization of monoidal categories, we would like to extend some of the properties desired of an ordinary monoidal product to the graded setting. We start by describing what it means for a graded monoidal product to be symmetric. Throughout the remainder of this section, we fix a graded monoidal category $(\C_\sbull, \mu, 1)$, using the same notation for the associators and unitors as above.

\begin{defn}
For each $m, n \in \NN$, we write $\tau_{m, n}$ for the \emph{transposition functor} $\C_m \times \C_n \to \C_n \times \C_m$. Explicitly, this is defined via the universal property of the product to be the unique functor making
\begin{equation*}
    \begin{tikzcd}[row sep = 2.5em]
        & \C_m \times \C_n \ar[dl] \ar[dr] \ar[d, dashed, "\tau_{m, n}" {description}]\\
        \C_n & \C_n \times \C_m \ar[l] \ar[r] & \C_m
    \end{tikzcd}
\end{equation*}
commute, where the unlabeled arrows are the appropriate projection functors.
\end{defn}

\begin{defn} \label{def:symmetric_monoidal}
We say $\C_\smallbullet$ is \emph{symmetric} if, for each $m, n \in \NN$, there exist natural isomorphisms $\sigma_{m, n} \colon \mu_{m, n} \cong \mu_{n, m} \tau_{m, n}$, called the \emph{braidings} of $\mu$, making the following diagrams commute for each $x, y, z \in \C$:
\begin{equation} \label{dgm:sym_unity}
    \begin{tikzcd}[row sep = 2.5em]
        & x\\
        \mu(1, x) \ar[ur, "{\lambda_x}"] \ar[rr, "{\sigma_{1, x}}"'] && \mu(x, 1) \ar[ul, "{\rho_x}"']
    \end{tikzcd}
\end{equation}
\begin{equation} \label{dgm:sym_assoc}
    \begin{tikzcd}[column sep = 4em, row sep = 2.3em]
        \mu(x, \mu(y, z)) \ar[r, "{\mu(x, \sigma_{y, z})}"] \ar[d, "{\alpha_{x, y, z}}"'] & \mu(x, \mu(z, y) \ar[d, "{\alpha_{x, z, y}}"]\\
        \mu(\mu(x, y), z) \ar[d, "{\sigma_{\mu(x, y), z}}"'] & \mu(\mu(x, z), y) \ar[d, "{\mu(\sigma_{x, z}, y)}"]\\
        \mu(z, \mu(x, y)) \ar[r, "{\alpha_{z, x, y}}"'] & \mu(\mu(z, x), y)
    \end{tikzcd}
\end{equation}
\begin{equation} \label{dgm:sym_involution}
    \begin{tikzcd}[row sep = 2.5em, column sep = 2.3em]
        \mu(x, y) \ar[dr, "{\sigma_{x, y}}"'] \ar[rr, equals] && \mu(x, y)\\
        & \mu(y, x) \period \ar[ur, "{\sigma_{y, x}}"']
    \end{tikzcd}
\end{equation}
\end{defn}

As in \Cref{def:graded_monoidal}, this definition directly generalizes the classical case. The coherence conditions for graded symmetric monoidal categories are precisely the same as the usual ones, with the exception that the sufficiently motivated reader could insert degree indices in each diagram above. We can likewise stratify the notion of a closed monoidal category by grading the internal hom.

\begin{defn} \label{def:closed_monoidal}
We say $\C_\sbull$ is \emph{closed} if, for each $m, n \in \NN$, there exists a bifunctor $\scF_{m, n} \colon \C_n^\op \times \C_{m + n} \to \C_m$ together with a natural isomorphism
\begin{equation*}
    \begin{tikzcd}[column sep = 4.5em, row sep = 3em]
        \C_m^\op \times \C_n^\op \times \C_{m + n} \ar[r, "\mu_{m, n} \times \C_{m + n}"] \ar[d, "\C_m^\op \times \scF_{m, n}"'] & \C_{m + n}^\op \times \C_{m + n} \ar[d, "{\C_{m + n}(-, -)}"] \ar[dl, Rightarrow, shorten >= 5ex, shorten <= 5ex, "\sim" {yshift = 1em, rotate = 27, xshift = -1em}] \\
        \C_m^\op \times \C_m \ar[r, "{\C_m(-, -)}"'] & \Set \period
    \end{tikzcd}
\end{equation*}
 The collection of bifunctors $\scF$ are called the \emph{internal hom} of $\C_\sbull$.
\end{defn}

\begin{rem}
In the case that $\mu$ is symmetric, $\scF$ defines a two-variable adjunction with $\mu$.
\end{rem}

\begin{rem} \label{rem:enriched_closed}
    If concerned with the $\V$-enriched setting, as outlined in \Cref{rem:enriched}, one simply replaces the hom-set functors $\C_k(-, -)\colon \C_k^\op \times \C_k \to \Set$ in \Cref{def:closed_monoidal} with the $\V$-enriched hom-functors $\underline{\C_k}(-, -) \colon \C_k^\op \times \C_k \to \V$.
\end{rem}

\subsection{Internalization}
\label{subsection:internalization}

In the preceding section, we laid the foundations for studying graded monoidal products. In this section, we describe the additional structure required to pull the external product of a graded monoidal category into the degree one component, and explore some of the basic behavior of this structure. Later, in \Cref{section:spectra}, we'll use these results to recover two constructions of point-set smash products of spectra, Day convolution and the smash product of $S$-modules, and show that these are closed and symmetric purely via the formal machinery of internalization.

\begin{defn} \label{def:internalization}
Let $(\C_\sbull, \mu, 1)$ be a graded monoidal category. An \emph{internalization of $\mu$} consists of:
\begin{enumerate}
    \item For each $m \in \NN$ a functor $\Phi_m \colon \C_m \to \C_1$, with a natural isomorphism $\iota$ which identifies $\Phi_1$ with the identity functor:
    \begin{equation*}
        \begin{tikzcd}[column sep = 3em]
            \C_1 \ar[r, bend right, equals, ""{name = B}] \ar[r, bend left, "\Phi_1"{name = A}] & \C_1 \period \ar[from = A, to = B, Rightarrow, shorten <=1ex, "\iota"]
        \end{tikzcd}
    \end{equation*} (As usual, we will write $\Phi$ to denote the collection of all the $\Phi_m$).
    \item For each $m, n \in \NN$, natural isomorphisms $\calL_{m, n}$ and $\calR_{m, n}$, called the \emph{left absorbers} and \emph{right absorbers} of $\Phi$, respectively\footnote{The absorbers are so-named because they allow the outermost factor of $\Phi$ to ``absorb" the inner ones.}:
    \begin{equation*}
        \begin{tikzcd}[column sep = 1em]
            & \C_1 \times \C_n \ar[rr, "\mu_{1, n}"] & {} \ar[d, Rightarrow,  shorten <= 1ex, shorten >= 0.5ex, "\calL_{m, n}"] & \C_{n + 1} \ar[dr, "\Phi_{n + 1}"]\\
            \C_m \times \C_n \ar[rr, "\mu_{m, n}"] \ar[dr, "\C_m \times \Phi_n"'] \ar[ur, "\Phi_m \times \C_n"] & {} & \C_{m + n} \ar[rr, "\Phi_{m + n}"] & {} & \C_1 \\
            & \C_m \times \C_1 \ar[rr, "\mu_{m, 1}"'] & {} \ar[u, Rightarrow, shorten <=1.8ex, shorten >= 0.4ex, "\calR_{m, n}"'] & \C_{m + 1} \period \ar[ur, "\Phi_{m + 1}"']
        \end{tikzcd}
    \end{equation*}
\end{enumerate} This data is further subject to the coherence conditions below -- given $x, y, z \in \C$, the following diagrams commute:

\begin{equation} \label{dgm:ambidexterity}
    \begin{tikzcd}[row sep = 2em]
        \Phi\mu(\Phi x, \Phi y) \ar[r, "\calL_{x, \Phi y}"] \ar[d, "\calR_{\Phi x, y}"'] & \Phi\mu(x, \Phi y) \ar[d, "\calR_{x, y}"]\\
        \Phi\mu(\Phi x, y) \ar[r, "\calL_{x, y}"'] & \Phi\mu(x, y),
    \end{tikzcd}
\end{equation}
\begin{equation} \label{dgm:right_assoc}
    \begin{tikzcd}[column sep = 3.5em]
        \Phi \mu(x, \Phi\mu(y, \Phi z)) \ar[r, "\calR_{x, \mu(y, \Phi z)}"] \ar[d, "{\Phi \mu(x, \calR_{y, z})}"'] & \Phi \mu(x, \mu(y, \Phi z)) \ar[r, "{\Phi \alpha_{x, y, \Phi z}}"] & \Phi \mu(\mu(x, y), \Phi z) \ar[d, "{\calR_{\mu(x, y), z}}"]\\
        \Phi \mu(x, \Phi \mu(y, z)) \ar[r, "{\calR_{x, \mu(y, z)}}"'] & \Phi \mu (x, \mu(y, z)) \ar[r, "{\Phi \alpha_{x, y, z}}"'] & \Phi \mu(\mu(x, y), z),
    \end{tikzcd}
\end{equation}
\begin{equation} \label{dgm:left_assoc}
    \begin{tikzcd}[column sep = 3.5em]
        \Phi\mu(\Phi \mu(\Phi x, y), z) \ar[r, "{\calL_{\mu(\Phi x, y), z}}"] \ar[d, "{\Phi \mu(\calL_{x, y}, z)}"'] & \Phi \mu(\mu(\Phi x, y), z) \ar[r, "{\Phi \alpha^{-1}_{\Phi x, y, z}}"] & \Phi \mu(\Phi x, \mu(y, z)) \ar[d, "{\calL_{x, \mu(y, z)}}"]\\
        \Phi\mu(\Phi \mu(x, y), z) \ar[r, "{\calL_{\mu(x, y), z}}"'] & \Phi \mu(\mu(x, y), z) \ar[r, "{\Phi\alpha_{x, y, z}}"'] & \Phi \mu(x, \mu(y, z)),
    \end{tikzcd}
\end{equation}
\begin{equation} \label{dgm:equilibrium}
    \begin{tikzcd}[column sep = 4em]
            \Phi \mu (x, \Phi \mu(y, z)) \ar[r, "{\calR_{x, \mu(y, z)}}"] \ar[d, "{\Phi \mu (x, \calL^{-1}_{y, z})}"'] & \Phi \mu(x, \mu(y, z)) \ar[r, "{\Phi \alpha_{x, y, z}}"] & \Phi \mu(\mu(x, y), z) \ar[d, "{\calL^{-1}_{\mu(x, y), z}}"]\\ 
            \Phi\mu(x, \Phi \mu(\Phi y, z)) \ar[d, "{\calR_{x, \mu(\Phi y, z)}}"'] && \Phi \mu(\Phi \mu(x, y), z) \ar[d, "{\Phi \mu(\calR^{-1}_{x, y}, z)}"]\\  
            \Phi \mu(x, \mu(\Phi y, z)) \ar[r, "{\Phi \alpha_{x, \Phi y, z}}"'] & \Phi \mu(\mu(x, \Phi y), z) \ar[r, "{\calL^{-1}_{\mu(x, \Phi y), z}}"'] & \Phi \mu(\Phi \mu(x, \Phi y), z),
    \end{tikzcd}
\end{equation} and lastly that
\begin{align}
    \calR_{m, 1} &= \Phi_{m + 1}\mu(\id, \iota) \label{eq:right_unity}\\
    \calL_{1, m} &= \Phi_{m + 1}\mu(\iota, \id). \label{eq:left_unity}
\end{align}

 If $(\C_\sbull, \mu, 1)$ is equipped with an internalization $\Phi$ as above, we say that $\C_\sbull$ is an \emph{internalized monoidal category} (we omit the term ``graded" for brevity -- occasionally we will refer to these as ``internalized categories"), and write $(\C_\sbull, \Phi)$ when $\Phi$ is not clear from context, leaving the absorbers implicit.
\end{defn}

As suggested by the name, the purpose of an internalization $\Phi$ on a graded monoidal category $\C_\smallbullet$ is to ``internalize" its external product $\mu$. The functors $\Phi$ allow us to normalize the degree of $\mu$, while the absorbers enable us to extract the ordinary monoidal coherence axioms from the graded ones.

\begin{thm}\label{thm:internalization}
Let $(\C_\sbull, \mu, 1)$ be a graded monoidal category equipped with an internalization $(\Phi, \calL, \calR)$. Then the bifunctor $\otimes \colon \C_1 \times \C_1 \to \C_1$ defined by
\begin{equation*}
    x \otimes y := \Phi_2 \mu_{1, 1}(x, y)
\end{equation*} defines a monoidal product on $\C_1$, with
\begin{enumerate}
    \item unit given by $I := \Phi_0(1),$
    \item associator $a$ given by the composite
    \begin{equation*}
        a_{x, y, z} := (\calL^{-1}_{2, 1})_{\mu(x, y), z} \circ \Phi_3 (\alpha_{1, 1, 1})_{x, y, z} \circ (\calR_{1, 2})_{x, \mu(y, z)},
    \end{equation*} for each $x, y, z \in \C_1$, depicted by the pasting diagram
    \begin{equation*}
        \begin{tikzcd}[column sep = 2.8em]
            & \C_1 \times \C_2 \ar[r, "\C_1 \times \Phi_2"] \ar[dr, "\mu_{1, 2}"] \ar[dd, Rightarrow, shorten <= 3.5ex, shorten >= 3.5ex, "\alpha_{1, 1, 1}"] & \C_1 \times \C_1 \ar[r, "\mu_{1, 1}"] & \C_2 \ar[dr, "\Phi_2"] \ar[dl, Rightarrow, shorten <= 3ex, shorten >= 3ex, "\calR_{1, 2}"] \\
            C_1 \times \C_1 \times \C_1 \ar[ur, "\C_1 \times \mu_{1, 1}"] \ar[dr, "\mu_{1, 1} \times \C_1"'] && \C_3 \ar[rr, "\Phi_3"'] & {} \ar[dl, Rightarrow, shorten <= 3ex, shorten >= 3ex, "\calL_{2,1}^{-1}"] & \C_1\\
            & \C_2 \times \C_1 \ar[r, "\Phi_2 \times \C_1"'] \ar[ur, "\mu_{2, 1}"'] &  \C_1 \times \C_1 \ar[r, "\mu_{1, 1}"'] & \C_2, \ar[ur, "\Phi_2"']
        \end{tikzcd}
    \end{equation*}
    \item and left and right unitors $\ell$ and $r$ (respectively), given by
    \begin{align*}
        \ell_x &:= \iota_x \circ \Phi_1(\lambda_1)_x \circ (\calL_{0, 1})_{1, x}\\
        r_x &:= \iota_x \circ \Phi_1(\rho_1)_x \circ (\calR_{1, 0})_{x, 1},
    \end{align*}
\end{enumerate} for each $x \in \C_1$, illustrated by the following pasting diagrams:

\adjustbox{scale=0.95,center}{%
    \begin{tikzcd}[column sep = 2.5em]
        & \C_0 \times \C_1 \ar[r, "\Phi_0 \times \C_1"{name = A'}] \ar[dr, "\mu_{0, 1}"] \ar[d, Rightarrow, shorten <= 1ex, shorten >= 1.5ex, "\lambda_1"] & \C_1 \times \C_1 \ar[r, "\mu_{1, 1}"] & \C_2 \ar[d, "\Phi_2"] \ar[dl, Rightarrow, shorten <= 3.5ex, shorten >= 3ex, "\calL_{0, 1}"]\\
        \C_1 \ar[rr, equals] \ar[ur, "{(1, \C_1)}"] \ar[rrr, equals, bend right, ""{name = B}] & {} & \C_1 \ar[Rightarrow, from = A', to = B, shorten <= 10ex, shorten >= 1ex, "\iota" {yshift = -1.8em}] \ar[r, "\Phi_1"'] & \C_1
    \end{tikzcd} \hspace{1em} 
    \begin{tikzcd}[column sep = 2.5em]
        \C_1 & \C_1 \ar[l, "\Phi_1"'] & {} & \C_1 \ar[ll, equals] \ar[lll, bend right, equals, ""{name = B}] \ar[dl, "{(\C_1, 1)}"]\\
        \C_2 \ar[u, "\Phi_2"] \ar[ur, Rightarrow, shorten <= 3.5ex, shorten >= 3ex, "\calR_{1, 0}"] & \C_1 \times \C_1 \ar[l, "\Phi_1"] & \C_1 \times \C_0 \period \ar[ul, "\mu_{1, 0}"] \ar[l, "\C_1 \times \Phi_0"{name = A}] \ar[u, Rightarrow, shorten <= 1ex, shorten >= 0.5ex, "\rho_1"] \ar[Rightarrow, from = A, to = B, shorten <= 10ex, shorten >= 1ex, "\iota"{yshift = 1.7em}]
    \end{tikzcd}
}

\noindent We say that $\otimes$ is the \emph{internalization of $\mu$ with respect to $\Phi$.}
\end{thm}
\begin{proof}
We must verify that the specified natural transformations satisfy the coherence conditions of a monoidal product, namely, that the following diagrams commute for all $x, y, z, w \in \C_1$:

\begin{equation}\label{dgm:monoidal_1}
    \begin{tikzcd}[column sep = 1.5em]
        & x \otimes y &\\
        x \otimes (1 \otimes y) \ar[ur, "{x \otimes \ell_y}"] \ar[rr, "{a_{x, 1, y}}"'] && (x \otimes 1) \otimes y \ar[ul, "{r_y \otimes y}"']
    \end{tikzcd} 
\end{equation}

\begin{equation}\label{dgm:monoidal_2}
    \begin{tikzcd}[column sep = 0em]
        x \otimes (y \otimes (z \otimes w)) \ar[rr, "{a_{x, y, z \otimes w}}"] \ar[d, "{x \otimes a_{y, z, w}}"'] && (x \otimes y) \otimes (z \otimes w) \ar[d, "{a_{x \otimes y, z, w}}"]\\
        x \otimes ((y \otimes z) \otimes w) \ar[dr, "{a_{x, y \otimes z, w}}"'] && ((x \otimes y) \otimes z) \otimes w \\
        & (x \otimes (y \otimes z)) \otimes w \period \ar[ur, "{a_{x, y, z} \otimes w}"']
    \end{tikzcd}
\end{equation}

\vspace{0.5em}

 First, we will show that \eqref{dgm:monoidal_1} commutes. By the definitions of $a$, $\ell$, and $r$, together with \eqref{eq:right_unity} and \eqref{eq:left_unity}, we can decompose the diagram into the outermost diagram below:

\adjustbox{scale=0.95,center}{%
    \begin{tikzcd}[column sep = 4em, row sep = 3em]
        & \Phi\mu(x,\Phi y) \ar[r, "{\calR_{x, y}}"] & \Phi \mu(x, y) & \Phi \mu(\Phi x, y) \ar[l, "{\calL_{x, y}}"'] \\
        \Phi \mu(x, \Phi \mu(1, y)) \ar[ur, "{\Phi\mu(x, \Phi \lambda_y)}"] \ar[r, "{\calR_{x, \mu(1, y)}}"'] & \Phi \mu(x, \mu(1, y)) \ar[rr, "{\Phi \alpha_{x, 1, y}}"'] \ar[ur, "{\Phi\mu(x, \lambda_y)}"'] && \Phi \mu(\mu(x, 1), y) \ar[ul, "{\Phi \mu(\rho_x, y)}"] & \Phi\mu(\Phi \mu(x, 1), y) \ar[l, "{\calL_{\mu(x, 1), y}}"] \ar[ul, "{\Phi\mu(\Phi \rho_x, y)}"']\\
        \Phi\mu(x, \Phi \mu(I, y)) \ar[u, "{\Phi\mu(x, \calL_{1, y})}"] \ar[r, "{\calR_{x, \mu(I, y)}}"'] & \Phi\mu(x, \mu(I, y)) \ar[rr, "{\Phi \alpha_{x, I, y}}"'] && \Phi\mu(\mu(x, I), y) \ar[r, "{\calL^{-1}_{\mu(x, I), y}}"'] & \Phi\mu(\Phi\mu(x, I), y) \period \ar[u, "{\Phi \mu(\calR_{x, 1}, y)}"']
    \end{tikzcd}
}

\vspace{0.5em}

\noindent But each inner subdiagram above commutes: the top left and right quadrilaterals commute by the naturality of $\calL$ and $\calR$; the top middle triangle is precisely the image of \eqref{dgm:unity} under $\Phi$; and the bottom rectangle is \eqref{dgm:equilibrium}. Therefore the entire diagram commutes.

It remains to be shown that \eqref{dgm:monoidal_2} commutes. Again, we decompose the diagram into subdiagrams to form the (somewhat menacing) one below:

\adjustbox{scale=0.72,center}{%
    \begin{tikzcd}[column sep = 5.2em, row sep = 3em]
        \Phi \mu(x, \Phi\mu(y, \Phi\mu(z, w))) \ar[r, "{\calR_{x, \mu(y, \Phi\mu(z, w))}}"] \ar[d, "{\Phi\mu(x, \calR_{y, \mu(z, w)})}"'] & \Phi \mu(x, \mu(y, \Phi \mu(z, w))) \ar[r, "{\Phi \alpha_{x, y, \Phi \mu(z, w)}}"] & \Phi \mu(\mu(x, y), \Phi \mu(z, w)) \ar[r, "{\calL^{-1}_{\mu(x, y), \Phi \mu(z, w)}}"] \ar[d, "{\calR_{\mu(x, y), \mu(z, w)}}"'] & \Phi \mu(\Phi \mu(x, y), \Phi \mu(z, w)) \ar[d, "{\calR_{\Phi\mu(x, y), \mu(z, w)}}"]\\
        \Phi \mu(x, \Phi\mu(y, \mu(z, w))) \ar[dd, "{\Phi \mu(x, \Phi \alpha_{y, z, w})}"'] \ar[r, "{\calR_{x, \mu(y, \mu(z, w))}}"'] & \Phi \mu(x, \mu(y, \mu(z, w))) \ar[dd, "{\Phi \mu(x, \alpha_{y, z, w})}"'] \ar[r, "{\Phi \alpha_{x, y, \mu(z, w)}}"'] & \Phi \mu(\mu(x, y), \mu(z, w)) \ar[r, "{\calL^{-1}_{\mu(x, y), \mu(z, w)}}"'] \ar[d, "{\Phi \alpha_{\mu(x, y), z, w}}"'] & \Phi \mu(\Phi \mu(x, y), \mu(z, w)) \ar[r, "{\Phi \alpha_{\Phi\mu(x, y), z, w}}"'] & \Phi \mu(\mu(\Phi\mu(x, y), z), w) \ar[d, "{\calL^{-1}_{\mu(\Phi\mu(x, y), z), w}}"]\\
        && \Phi\mu(\mu(\mu(x, y), z), w) \ar[r, "{\calL^{-1}_{\mu(\mu(x, y), z), w}}"'] & \Phi \mu(\Phi \mu(\mu(x, y), z), w) \ar[r, "{\Phi\mu(\calL^{-1}_{\mu(x, y), z}, w)}"'] & \Phi \mu(\Phi\mu(\Phi\mu(x, y), z), w)\\
        \Phi \mu(x, \Phi \mu(\mu(y, z), w)) \ar[r, "{\calR_{x, \mu(\mu(y, z), w)}}"'] \ar[d, "{\Phi\mu(x, \calL^{-1}_{\mu(y, z), w})}"'] & \Phi\mu(x, \mu(\mu(y, z), w)) \ar[r, "{\Phi \alpha_{x, \mu(y, z), w}}"'] & \Phi \mu(\mu(x, \mu(y, z)), w) \ar[u, "{\Phi\mu(\alpha_{x, y, z}, w)}"] \ar[r, "{\calL^{-1}_{\mu(x, \mu(y, z)), w}}"'] & \Phi\mu(x, \mu(\Phi \mu(y, z), w)) \ar[r, "{\Phi\mu(\Phi \alpha_{x, y, z}, w)}"'] & \Phi\mu(\Phi \mu(\mu(x, y), z), w) \ar[u, "{\Phi \mu(\calL^{-1}_{\mu(x, y), z}, w)}"']\\
        \Phi \mu(x, \Phi \mu(\Phi \mu(y, z), w)) \ar[r, "{\calR_{x, \mu(\Phi \mu(y, z), w)}}"'] & \Phi\mu(x, \mu(\Phi \mu(y, z), w)) \ar[r, "{\Phi \alpha_{x, \Phi \mu(y, z), w}}"'] & \Phi \mu(\mu(x, \Phi \mu(y, z)), w) \ar[r, "{\calL^{-1}_{\mu(x, \Phi\mu(y, z)), w}}"'] & \Phi \mu(\Phi \mu(x, \Phi \mu(y, z)), w) \period \ar[u, "{\Phi\mu(\calR_{x, \mu(y, z)}, w)}"']
    \end{tikzcd}
}

\vspace{0.5em}

\noindent Once again, we will find that each subdiagram commutes. The top-left hexagon is just \eqref{dgm:right_assoc}, the two middle-right hexagons are \eqref{dgm:left_assoc}, and the top-right square is \eqref{dgm:ambidexterity}. Additionally, the middle-left square commutes by the naturality of $\calR$. Finally, the central pentagon is the image of \eqref{dgm:associativity} under $\Phi$, and thus also commutes. Therefore the entire diagram commutes, completing the proof.
\end{proof}

By imposing some additional structure on an internalization, we can also recover standard properties of the internalized monoidal product like symmetry and closure.

\begin{defn} \label{def:symmetric_internalization}
Let $(\C_\sbull, \mu, 1)$ be a symmetric graded monoidal category with braiding $\sigma$, equipped with an internalization $(\Phi, \calL, \calR)$. Then $\Phi$ is a \emph{symmetric internalization} if, for all $x, y \in \C$, the diagram
\begin{equation*}
    \begin{tikzcd}[row sep = 2.5em]
        \Phi\mu(\Phi x, y) \ar[r, "{\Phi \sigma_{\Phi x, y}}"] \ar[d, "{\calL_{x, y}}"'] & \Phi\mu(y, \Phi x) \ar[d, "{\calR_{y, x}}"]\\
        \Phi\mu(x, y) \ar[r, "{\Phi \sigma_{x, y}}"'] & \Phi\mu(y, x)
    \end{tikzcd}
\end{equation*} commutes.
\end{defn}

\begin{cor} \label{cor:symmetric_internalization}
Let $(\C_\sbull, \mu, 1)$ be a symmetric graded monoidal category with braiding $\sigma$, equipped with a symmetric internalization $(\Phi, \calL, \calR)$. Then the natural isomorphism
\begin{equation*}
    s_{x, y} := \Phi_2 (\sigma_{1, 1})_{x, y}
\end{equation*} defines a braiding on the internalization $\otimes$ of $\mu$ with respect to $\Phi$, making $(\C_1, \otimes, I)$ into a symmetric monoidal category.
\end{cor}
\begin{proof}
We must verify that the natural isomorphism $s$ satisfies the three coherence conditions for a symmetric monoidal category, namely, that the following diagrams commute for all $x, y, z \in \C_1$:

\begin{equation}\label{dgm:braiding_1}
    \begin{tikzcd}[row sep = 2.5em]
        & x &\\
        I \otimes x \ar[ur, "\ell_x"] \ar[rr, "s_{1, x}"'] && x \otimes I \ar[ul, "r_x"']
    \end{tikzcd}
\end{equation}

\begin{equation} \label{dgm:braiding_2}
    \begin{tikzcd}[row sep = 2.5em, column sep = 3.5em]
        x \otimes (y \otimes z) \ar[r, "{x \otimes s_{y, z}}"] \ar[d, "{a_{x, y, z}}"'] & x \otimes (z \otimes y) \ar[d, "{a_{x, z, y}}"]\\
        (x \otimes y) \otimes z \ar[d, "{s_{x \otimes y, z}}"'] & (x \otimes z) \otimes y \ar[d, "{s_{x, z} \otimes y}"]\\
        z \otimes (x \otimes y) \ar[r, "{a_{z, x, y}}"'] & (z \otimes x) \otimes y
    \end{tikzcd}
\end{equation}

\begin{equation} \label{dgm:braiding_3}
    \begin{tikzcd}[row sep = 2.5em]
        x \otimes y \ar[dr, "{s_{x, y}}"'] \ar[rr, equals] && x \otimes y\\
        & y \otimes x, \ar[ur, "{s_{y, x}}"']
    \end{tikzcd}
\end{equation} where $a, \ell, $ and $r$ are defined in \Cref{thm:internalization}. First, we observe that \eqref{dgm:braiding_1} commutes by the commutativity of
\begin{equation*}
    \begin{tikzcd}[row sep = 2.5em, column sep = 1.5em]
        & \Phi x \\
        \Phi \mu(1, x) \ar[ur, "{\Phi \lambda_x}"] \ar[rr, "{\Phi \sigma_{1, x}}"] && \Phi \mu(x, 1) \ar[ul, "{\Phi \rho_x}"']\\
        \Phi \mu(\Phi (1), x) \ar[u, "{\calL_{1, x}}"] \ar[rr, "{\Phi \sigma_{\Phi(1), x}}"'] && \Phi \mu(x, \Phi(1)) \period \ar[u, "{\calR_{x, 1}}"']
    \end{tikzcd}
\end{equation*}
 
\noindent Each subdiagram commutes---the triangle is image of \eqref{dgm:sym_unity} under $\Phi$, and the square commutes by \Cref{def:symmetric_internalization}---and post-composition with $\iota_x$ yields \eqref{dgm:braiding_1}.

Next, we will show that \eqref{dgm:braiding_2} commutes. As usual, we begin by decomposing the diagram into smaller subdiagrams:

\adjustbox{scale=1,center}{%
    \begin{tikzcd}[column sep = 4.5em]
        & \Phi \mu(x, \Phi\mu(y, z)) \ar[r, "{\Phi \mu(x, \Phi \sigma_{y, z})}"] \ar[d, "{\calR_{x, \mu(y, z)}}"'] & \Phi \mu(x, \Phi \mu(z, y)) \ar[d, "{\calR_{x, \mu(z, y)}}"]\\
        & \Phi \mu(x, \mu(y, z)) \ar[r, "{\Phi \mu(x, \sigma_{y, z})}"'] \ar[d, "{\Phi \alpha_{x, y, z}}"'] & \Phi \mu(x, \mu(z, y)) \ar[d, "{\Phi\alpha_{x, z, y}}"]\\
        \Phi\mu(\Phi\mu(x, y), z) \ar[d, "{\Phi\sigma_{\Phi\mu(x, y), z}}"'] & \Phi \mu(\mu(x, y), z) \ar[l, "{\calL^{-1}_{\mu(x, y), z}}"'] \ar[d, "{\Phi \sigma_{\mu(x, y), z}}"'] & \Phi \mu(\mu(x, z), y) \ar[r, "{\calL^{-1}_{\mu(x, z), y}}"] \ar[d, "{\Phi \mu(\sigma_{x, z}, y)}"] & \Phi \mu(\Phi \mu(x, z), y) \ar[d, "{\Phi \mu(\Phi\sigma_{x, z}, y)}"]\\
        \Phi\mu(z, \Phi\mu(x, y)) \ar[r, "{\calR_{z, \mu(x, y)}}"'] & \Phi\mu(z, \mu(x, y)) \ar[r, "{\Phi\alpha_{z, x, y}}"'] & \Phi \mu(\mu(z, x), y) \ar[r, "{\calL^{-1}_{\mu(z, x), y}}"'] & \Phi\mu(\Phi\mu(z, x), y) \period
    \end{tikzcd}
}

\vspace{0.5em}

\noindent The topmost square commutes by the naturality of $\calR$. The bottom-left square commutes by \Cref{def:symmetric_internalization}. The bottom-right square commutes by the naturality of $\calL$. Finally, the middle hexagon commutes since it is the image of \eqref{dgm:sym_assoc} under $\Phi$. Thus, \eqref{dgm:braiding_2} commutes.

Finally, it is immediate that \eqref{dgm:braiding_3} commutes, since it is just the image of \eqref{dgm:sym_involution} under $\Phi$. Therefore $s$ defines a braiding for $\otimes$ which makes $(\C_1, \otimes, I)$ into a symmetric monoidal category.
\end{proof}

\begin{defn} \label{def:closed_internalization}
Let $(\C_\sbull, \mu, 1)$ be a closed graded monoidal category with internal hom $\scF$, equipped with an internalization $(\Phi, \calL, \calR)$. Then $\Phi$ is a \emph{closed internalization} if, for each $m \in \NN$, there exists a right adjoint $\Psi_m \colon \C_1 \to \C_m$ to $\Phi_m$. For brevity, we write $\Phi \dashv \Psi$ to denote the adjunction.
\end{defn}

\begin{cor} \label{cor:closed_internalization}
Let $(\C_\sbull, \mu, 1)$ be a closed graded monoidal category with internal hom $\scF$, equipped with a closed internalization $(\Phi, \calL, \calR)$, with $\Phi \dashv \Psi$. Then $(\C_1, \otimes, I)$ is a closed monoidal category with internal hom
\begin{equation*}
    F(x, y) := \scF_{1, 1}(x, \Psi_2 y)
\end{equation*} for each $x, y \in \C_1$.
\end{cor}
\begin{proof}
Let $x, y, z \in \C_1$. Then the claim follows immediately by composing the natural isomorphisms
\begin{align*}
     \C_1(x \otimes y, z) &= \C_1(\Phi_2\mu_{1, 1}(x, y), z)\\
     &\cong \C_2(\mu_{1, 1}(x, y), \Psi_2 z)\\
     &\cong \C_1(x, \scF_{1, 1}(y, \Psi_2 z))\\
     &= \C_1(x, F(y, z)).
\end{align*} Diagramatically, this is expressed by the pasting diagram below (where, as noted in \Cref{rem:enriched_closed}, one can replace $\Set$ with the appropriate base category in the enriched setting):
\begin{equation*}
    \begin{tikzcd}[column sep = 3em, /tikz/column 3/.style={column sep=4em}, /tikz/column 1/.style={column sep=1em}]
        & \C_2^\op \times \C_1 \ar[r, "\Phi_2 \times \C_1"] \ar[dr, "\C_2^\op \times \psi_2"'{yshift = 0.7ex}] \ar[dd, equals, shorten <= 5ex, shorten >= 5ex] & \C_1^\op \times \C_1 \ar[dr, "{\C_1(-, -)}"] \ar[d, Rightarrow, shorten <= 0.5ex, shorten >= 1ex, "\sim" {rotate = 90, yshift = -1ex, xshift = -0.8ex}] \\
        \C_1^\op \times \C_1^\op \times \C_1 \ar[ur, "\mu_{1, 1} \times \C_1"] \ar[dr, "\C_1^\op \times \C_1^\op \times \psi_2"'] && \C_2^\op \times \C_2 \ar[r, "{\C_2(-, -)}"{xshift = -1ex}] \ar[d, Rightarrow, shorten <= 0.5ex, shorten >= 1ex, "\sim" {rotate = 90, yshift = -1ex, xshift = -0.8ex}] & \Set\\
        & \C_1^\op \times \C_1^\op \times \C_2 \ar[ur, "\mu_{1, 1} \times \C_2"{yshift = -0.5ex}] \ar[r, "\C_1^\op \times \scF_{1,1}"'] & \C_1^\op \times \C_1 \period \ar[ur, "{\C_1(-, -)}"']
    \end{tikzcd}
\end{equation*}
\end{proof}

\section{Categorical Properties of Internalization}
\label{section:categorical_properties}

In this section, we outline the basic (higher) categorical properties of graded monoidal categories and internalized categories. With suitable notions of functors and natural transformations, we will see that these form 2-categories equipped with compatible adjunctions with classical structures. In particular, internalization transforms a graded monoidal category \emph{functorially} into an ordinary monoidal category (\Cref{prop:int-functorial}), and every monoidal category can be lifted to a \emph{free} internalized category whose internalization recovers the original monoidal product in degree one (\Cref{thm:int-section}, \Cref{thm:free-forget-adj-mon-int}).

\subsection{The 2-Category of Graded Monoidal Categories}

We begin by defining (lax, strong, and strict) graded monoidal functors, which is directly analogous to the corresponding definition(s) for ordinary monoidal categories \cite[Chapter~XI]{maclane}.

\begin{defn} \label{def:gr-mon-functor}
A \emph{lax graded monoidal functor} $f_\sbull \colon \C_\sbull \to \D_\sbull$ consists of
\begin{enumerate}
    \item a functor $f_m \colon \C_m \to \D_m$ for each $m \in \mathbb{N}$,
    \item a natural transformation $$\zeta^{m,n} \colon \mu^\D_{m, n}(f_m, f_n) \Rightarrow f_{m+n}\mu^\C_{m,n}$$ for each $m, n \in \mathbb{N}$,
    \item a morphism $\nu \colon 1^\D \to f_0(1^\C)$,
\end{enumerate}such that the following diagrams commute (for all $x \in \C_m, y \in \C_n,$ and $ z \in \C_r$):
\begin{equation*}
\begin{tikzcd}[column sep = 2.5em, row sep = 2.5em]
    \mu^\D(f_m(x), \mu^\D(f_n(y), f_r(z))) \ar[r, "\alpha^\D"] \ar[d, "\mu^\D({\id, \zeta^{n, r}})"'] & \mu^\D(\mu^\D(f_m(x), f_n(y)), f_r(z)) \ar[d, "\mu^\D({\zeta^{m, n}, \id})"]\\
    \mu^\D(f_m(x), f_{n+r}\mu^\C(y, z)) \ar[d, "\zeta^{m, n+r}"'] & \mu^\D(f_{m+n}\mu^\C(x, y), f_r(z)) \ar[d, "\zeta^{m+n, r}"] \\
    f_{m+n+r}\mu^\C(x, \mu^\C(y, z))\ar[r, "f_{m+n+r}\alpha^\C"'] & f_{m+n+r}\mu^\C(\mu^\C(x, y), z)
\end{tikzcd}
\end{equation*}

\begin{equation*}
\adjustbox{scale=1.0,center}{%
\begin{tikzcd}
    \mu^\D(1^\D, f_m(x)) \ar[r, "\mu^\D({\nu,} \id)"] \ar[d, "\lambda^\D"'] & \mu^\D(f_0(1^\C), f_m(x)) \ar[d, "\zeta^{0, m}"] \\
    f_m(x) & \ar[l, "f_m\lambda^\C"] f_m\mu^\C(1^\C, x)
\end{tikzcd} \hspace{2em}
\begin{tikzcd}
    \mu^\D(f_m(x), 1^\D) \ar[r, "\mu^\D(\id{,} \nu)"] \ar[d, "\rho^\D"'] & \mu^\D(f_m(x), f_0(1^\C)) \ar[d, "\zeta^{m, 0}"] \\
    f_m(x) & \ar[l, "f_m\rho^\C"] f_m\mu^\C(x, 1^\C) \period
\end{tikzcd}
}
\end{equation*}

If $\nu$ and each $\zeta^{m,n}$ are isomorphisms then $f_\sbull$ is a \emph{strong graded monoidal functor}, and if they are identity morphisms then $f_\sbull$ is a \emph{strict graded monoidal functor}.
\end{defn}

The identity $\mathrm{id}_\sbull \colon \C_\sbull \to \C_\sbull$ is defined to be the componentwise identity in each degree, and is a strict graded monoidal functor. The composite $(gf)_\sbull$ of a pair of lax (resp. strong, strict) graded monoidal functors $f_\sbull \colon \C_\sbull \to \D_\sbull$, $g_\sbull \colon \D_\sbull \to \E_\sbull$ is obtained by componentwise composition over the grading, while the composites
$$
\mu^\E_{m,n}(g_mf_m(x), g_nf_n(y)) \to g_{m+n}(\mu^\D_{m,n}(f_m(x), f_n(y))) \to g_{m+n}f_{m+n}\mu_{m,n}^\C(x, y)
$$
and
$$
1^\E \to g_0(1^\D) \to g_0f_0(1^\C)
$$
define the structure maps making $(gf)_\sbull$ a lax (resp. strong, strict) functor. We elide the verification of the coherence diagrams; the argument is identical to the classical case, except extra care must be taken to track the degrees of each object.

\begin{rem}
While the choice to favor any one strength of graded monoidal functor is a matter of taste, for our purposes it will be largely irrelevant. We will henceforth assume all (graded) monoidal functors are lax, but note that throughout this section they could be made strong or strict without changing any of the subsequent results. In light of this remark, and for brevity, we will simply refer to these as ``graded monoidal functors."
\end{rem}

\begin{defn}
Graded monoidal categories\footnote{As usual when referring to categories of categories, we must fix an upper bound on the size of the categories in question to avoid size issues. To be precise, $\GrMonCat$ is the category of \emph{small} graded monoidal categories (wherein each component is a small category).} and graded monoidal functors between them form a category, which we denote $\GrMonCat$.
\end{defn}

As is typical for categories whose morphisms are some flavor of functor, a more complete description of $\GrMonCat$ is as a 2-category \cite[Section XII.3]{maclane}. To make sense of this, we first need a notion of natural transformations between graded monoidal functors, which, predictably, is just the graded version of monoidal transformations \cite[Chapter XI]{maclane}.

\begin{defn}
Given graded monoidal functors $f_\sbull, g_\sbull \colon C_\sbull \to D_\sbull$, a \emph{graded monoidal (natural) transformation} $\beta_\sbull \colon f_\sbull \Rightarrow g_\sbull$ consists of a natural transformation $\beta_m \colon f_m \Rightarrow g_m$ for each $m \in \mathbb{N}$ such that the following diagrams commute (for all $x \in \C_m, y \in \C_n$):
\begin{equation*}
\begin{tikzcd}
    & 1^\D \ar[dl, "\nu"'] \ar[dr, "\nu"] \\
    f_0(1^\C) \ar[rr, "\beta_0"'] && g_0(1^\C)
\end{tikzcd}
\end{equation*}
\begin{equation*}
\begin{tikzcd}[column sep = 4em, row sep = 2.5em]
    \mu^\D(f_m(x), f_n(y)) \ar[r, "\mu^\D(\beta_m{,} \beta_n)"] \ar[d, "\zeta^{m, n}"'] & \mu^\D(g_m(x), g_n(y)) \ar[d, "\zeta^{m, n}"] \\
    f_{m+n}\mu^\C(x, y) \ar[r, "\beta_{m+n}"'] & g_{m+n}\mu^\C(x, y)
\end{tikzcd}
\end{equation*}
We abuse notation slightly by denoting the structure maps of both $f_\sbull$ and $g_\sbull$ using the same symbols $\nu$ and $\zeta$, in an admittedly futile attempt at minimizing clutter. Throughout this section, when it is clear which structure maps are being referenced, we will renew this abuse.
\end{defn}

With graded monoidal functors as 1-morphisms and graded monoidal natural transformations as 2-morphisms, $\GrMonCat$ forms a 2-category. Since graded monoidal functors and natural transformations are defined componentwise, isolating any single degree $m$ defines a 2-functor $\GrMonCat \to \Cat,$ where $\Cat$ denotes the 2-category of (small) categories. We give special attention to the case when $m = 1$:

\begin{prop} \label{prop:forget-gr}
There is a 2-functor $\Forget_\Gr \colon \GrMonCat \to \Cat$ which sends a graded monoidal category $(\C_\sbull, \mu, 1)$ to its degree one component $\C_1$.
\end{prop}

Much like one can construct a free monoidal category from any category via the finite list monad \cite[Section VII.3, Theorem 2]{maclane}, there is an analogous \emph{free graded monoidal product} on any category, where the grading is given by list length. As one might expect, this construction is 2-functorial.

\begin{prop} \label{prop:free-gr}
There is a 2-functor $\Free_\Gr \colon \Cat \to \GrMonCat$ which carries a category $\C$ to the graded monoidal category $\Free_\Gr(\C)_\sbull$ defined by
\begin{align*}
    \Free_\Gr(\C)_m &:= \C^m
\end{align*} equipped with the empty string $1^\Free \in \C^0$ as the unit $1^\Free$ and graded monoidal product $\mu^\Free$ given by concatenation:
\begin{equation*}
    \mu^\Free_{m, n}((x_1, \dots, x_m), (y_1, \dots, y_n)) := (x_1, \dots, x_m, y_1, \dots, y_n).
\end{equation*}
\end{prop}

\begin{proof}
The structure maps for the graded monoidal product $\mu^\Free$ are trivial (that is, $\Free_\Gr(\C)_\sbull$ is a \emph{strict} graded monoidal product), and thus the axioms are free.

For a functor $f \colon \C \to \D$ we define 
$$
\Free_\Gr(f)_m := f^m \colon \C^m \to \D^m
$$that is, by termwise application of $f$, which is clearly compatible with the concatenation product and empty string unit of $\Free_\Gr(\C)_\sbull$. 

Similarly, for a natural transformation $\beta \colon f \Rightarrow g$, let $\Free_\Gr(\beta)_m := \beta^m$, which one easily checks defines a graded monoidal transformation $\Free_\Gr(f)_\sbull \Rightarrow \Free_\Gr(g)_\sbull$. Verifying that $\Free$ is 2-functorial is straightforward.
\end{proof}

This highlights the significance of the degree one component -- lists of length one contain exactly the data of the original category, which is why we view $\Forget_\Gr$ as a kind of ``forgetful" functor.

\begin{prop} \label{prop:free-forget-gr-section}
    The functor $\Free_\Gr$ is a section of $\Forget_\Gr.$
\end{prop}
\begin{proof}
From the definitions in Propositions \ref{prop:forget-gr} and \ref{prop:free-gr} we verify directly that
\begin{equation*}
    (\Forget_\Gr \circ \Free_\Gr )(\C) = \Free_\Gr(\C)_1 = \C^1 = \C.
\end{equation*} on objects $\C \in \Cat$. The cases for 1-morphisms (functors) and 2-morphisms (natural transformations) are identical.
\end{proof}

The free graded monoidal category above is indeed ``free," in the sense that it is left adjoint to $\Forget_\Gr$.

\begin{prop} \label{prop:free-forget-adj-cat-grmoncat}
There is a (strict) 2-adjunction of 2-functors 
\begin{equation*}
    \begin{tikzcd}
        \Cat \ar[rr, bend left = 20, "\Free_\Gr" {xshift=1.5em}] \ar[rr, phantom, "\bot"{xshift = 0.6em}] && \GrMonCat \period \ar[ll, bend left = 20, "\Forget_\Gr" {xshift=0.5em}]
    \end{tikzcd}
\end{equation*}
\end{prop}

\begin{proof}
To establish a 2-adjunction $\Free_\Gr \dashv \Forget_\Gr,$ we construct the unit $\eta$ and counit $\epsilon$ as follows: in light of Proposition \ref{prop:free-forget-gr-section}, we define the component $\eta_\C$ to simply be the identity functor on $\C$. Given $(\D_\sbull, \mu^\D, 1^\D) \in \GrMonCat$, we define the $m\Th$ graded component of the counit $(\epsilon_{\D_\sbull})_m \colon \D_1^m \to \D_m$ by
$$
(\epsilon_{\D_\sbull})_m(x_1, \dots, x_m) := \mu^\D(x_1, \dots, x_m)
$$(where we fix a bias for the $m$-fold iterated $\mu^\D$ product, say by associating to the right) if $m > 1$ and take
\begin{align*}
(\epsilon_{\D_\sbull})_1(x) &:= x\\
(\epsilon_{\D_\sbull})_0(\ast) &:= 1^\D.
\end{align*}Then for $(x_1, \dots, x_m) \in \C^m = \Free_\Gr(\C)_m$ we have
\begin{align*}
(\epsilon_{\Free_\Gr(\C)_\sbull} \circ \Free_\Gr(\eta_\C))_m(x_1, \dots, x_m) &= \epsilon_{\Free_\Gr(\C)_\sbull}(\eta_\C(x_1), \dots, \eta_\C(x_m))\\
&= \mu^\Free(x_1, \dots, x_m)\\
&= (x_1, \dots, x_m),
\end{align*}so $$\epsilon_{\Free_\Gr(\C)_\sbull} \circ \Free_\Gr(\eta_\C)_\sbull = \id_{\Free_\Gr(\C)_\sbull}.$$ On the other hand, for $y \in \D_1 = \Forget_\Gr(\D_\sbull)$ we have
$$
(\Forget_\Gr(\epsilon_{\D_\sbull}) \circ \eta_{\Forget_\Gr(\D_\sbull)})(y) = \Forget_\Gr(\epsilon_{\mathrm{\D_\sbull}})(y) = ({\epsilon_{\D_\sbull}})_1(y) = y.
$$An identical argument verifies the equations for morphisms, so indeed $\eta$ and $\epsilon$ define a (strict) 2-adjunction $\Free_\Gr \dashv \Forget_\Gr.$ 
\end{proof}

\begin{cor}\label{cor:gr-coreflective-cat}
    The functor $\Free_\Gr$ embeds $\Cat$ as a coreflective subcategory of $\GrMonCat$.
\end{cor}
\begin{proof}
    The functor $\Free_\Gr$ is injective on objects, since we can recover the category $\C$ from the degree one component of $\Free_\Gr(\C)_\sbull$; the same is clearly true for functors and natural transformations. Thus $\Free_\Gr$ is an embedding of 2-categories. The unit of the adjunction in Proposition \ref{prop:free-forget-adj-cat-grmoncat} is the identity (so in particular is an isomorphism), making this embedding coreflective.
\end{proof}

\subsection{The 2-Category of Internalized Categories}

Having established the relationship between the 2-categories $\Cat$ and $\GrMonCat$, we now turn our attention to \emph{internalized} graded monoidal categories. As one might hope, these too assemble into a 2-category $\IntMonCat$. This requires an appropriate extension of graded monoidal functors and natural transformations to accomodate the additional structure of an internalization.

\begin{defn}
Let $\C_\sbull$ and $\D_\sbull$ be graded monoidal categories equipped with internalizations $\Phi^\C$ and $\Phi^\D$, respectively. A \emph{lax internalized (monoidal) functor} $f_\sbull \colon (\C_\sbull, \Phi^\C) \to (\D_\sbull, \Phi^\D)$ is a lax graded monoidal functor $\C_\sbull \to \D_\sbull$ together with natural transformations
$$
\xi^m \colon \Phi^\D_m f_m \Rightarrow f_1\Phi^\C_m
$$such that the following diagrams commute (for all $x \in \C_m, y \in \C_n,$ and $z \in \C_r$):
\begin{equation}
\begin{tikzcd}[column sep=-1em]
    \Phi_{n+1}^\D\mu^\D(\Phi_m^\D f_m(x), f_n(y)) \ar[rr, "{\Phi_{n+1}^\D\mu^\D(\xi^m, \id)}"] \ar[d, "{\calL^\D}"'] && \Phi_{n+1}^\D\mu^\D(f_1\Phi_m^\C(x), f_n(y)) \ar[d, "\Phi_{n+1}^\D\zeta^{1,n}"]\\
    \Phi^\D_{m+n}\mu^\D(f_m(x), f_n(y)) \ar[d, "{\Phi^\D_{m+n}\zeta^{m,n}}"'] && \Phi_{n+1}^\D f_{n+1} \mu^\C(\Phi_m^\C(x), y) \ar[d, "\xi^{n+1}"]\\
    \Phi^\D_{m+n}f_{m+n}\mu^\C(x,y) \ar[dr, "\xi^{m+n}"'] && f_1\Phi_{n+1}^\C \mu^\C(\Phi_m^\C(x), y) \ar[dl, "f_1\calL^\C"]\\
    & f_1\Phi^\C_{m+n}\mu^\C(x,y)
\end{tikzcd} \label{diag:lax-coh-left}
\end{equation}

\begin{equation}
\begin{tikzcd}[column sep=-1em]
    \Phi_{m+1}^\D\mu^\D(f_m(x), \Phi_n^\D f_n(y)) \ar[rr, "{\Phi_{m+1}^\D\mu^\D(\id, \xi^n)}"] \ar[d, "{\calR^\D}"'] && \Phi_{m+1}^\D\mu^\D(f_m(x), f_1\Phi_n^\C(y)) \ar[d, "\Phi_{m+1}^\D\zeta^{m,1}"]\\
    \Phi^\D_{m+n}\mu^\D(f_m(x), f_n(y)) \ar[d, "{\Phi^\D_{m+n}\zeta^{m,n}}"'] && \Phi_{n+1}^\D f_{n+1} \mu^\C(\Phi_m^\C(x), y) \ar[d, "\xi^{n+1}"]\\
    \Phi^\D_{m+n}f_{m+n}\mu^\C(x,y) \ar[dr, "\xi^{m+n}"'] && f_1\Phi_{m+1}^\C \mu^\C(x, \Phi_n^\C(y)) \ar[dl, "f_1\calR^\C"]\\
    & f_1\Phi^\C_{m+n}\mu^\C(x,y)
\end{tikzcd} \label{diag:lax-coh-right}
\end{equation}

\begin{equation}
\begin{tikzcd}[column sep=1em]
    & \Phi^\D_1 f_1 \Phi^\C_1(x) \ar[dl, "\Phi^\D_1 f_1 \iota^\C"'] \ar[dr, "\iota^\D"] \\
    \Phi^\D_1 f_1(x) \ar[rr, "\xi^1"'] && f_1(x) \period
\end{tikzcd} \label{diag:lax-coh-triv}
\end{equation}

If $f_\sbull$ is a strong graded monoidal functor and each natural transformation $\xi^m$ is an isomorphism we say $f_\sbull$ is a \emph{strong internalized (monoidal) functor}; if $f_\sbull$ is a strict graded monoidal functor and the $\xi^m$ are all identities then $f_\sbull$ is a \emph{strict internalized (monoidal) functor}.
\end{defn}

Going forward, we assume the same choice of strength for all internalized functors as graded monoidal functors. As above, we default to the weakest case of all lax functors, but note that if every monoidal functor in sight is strengthened to the strong or strict cases the following results remain valid.

Given an internalized category $\C_\sbull$, we denote by $\Int(\C_\sbull)$ the corresponding monoidal category equipped with the internalized product of \Cref{thm:internalization}. This construction extends to internalized functors:

\begin{prop} \label{prop:int_of_functor_is_monoidal}
    Let $f_\sbull \colon (\C_\sbull, \Phi^\C) \to (\D_\sbull, \Phi^\D)$ be an internalized functor. Then, equipping $\C_1$ (resp.\ $\D_1$) with the internalized monoidal product of \Cref{thm:internalization}, the functor $\Int(f) := f_1 \colon \C_1 \to \D_1$ is a monoidal functor.
\end{prop}
\begin{proof}
First, we must define structure maps
\begin{align*}
    \zeta^\Int &\colon \Phi^\D_2\mu^\D_{1,1}(f, f) \Rightarrow f\Phi^\C_2\mu_{1,1}\\
    \nu^\Int & \colon \Phi^\D_0(1^\D) \to f_1(\Phi^\C_0(1^\C))
\end{align*}for $\Int(f)$. Given $x, y \in \C_1$, we let $\zeta^\Int_{x,y}$ be the composite
\begin{equation*}
    \Phi^\D_2\mu^\D_{1,1}(f_1(x), f_1(y)) \xrightarrow{\Phi^\D\zeta^{1,1}_{x,y}} \Phi^\D_2f_2 \mu_{1,1}^\C(x,y) \xrightarrow{\xi^2_{\mu^\C(x,y)}} f_1\Phi^\C_2\mu_{1,1}^\C(x,y),
\end{equation*}and we take $v^\Int$ to be the morphism
\begin{equation*}
    \Phi^\D_0(1^\D) \xrightarrow{\Phi^\D_0\nu} \Phi^\D_0(f_0(1^\C)) \xrightarrow{\xi^0_{1^\C}} f_1\Phi_0^\C(1^\C).
\end{equation*} Next, we must check that these structure maps satisfy the coherence conditions of a monoidal functor, namely, that the following diagrams commute (for all $x, y, z \in \C_1$):
\begin{equation*}
\begin{tikzcd}
    \Phi^\D_2 \mu^\D(f(x), \Phi^\D_2\mu^\D(f(y), f(z))) \ar[r, "a^\D"] \ar[d, "{\Phi^\D_2 \mu^\D(\id, \zeta^\Int)}"'] & \Phi^\D_2\mu^\D(\Phi^\D_2\mu^\D(f(x), f(y)), f(z)) \ar[d, "{\Phi^\D_2 \mu^\D(\zeta^\Int, \id)}"] \\
    \Phi^\D_2 \mu^\D(f(x), f\Phi^\C_2\mu^\C(y, z)) \ar[d, "\zeta^\Int"'] & \Phi^\D_2\mu^\D(f\Phi^\C_2 \mu^\C(x, y), f(z)) \ar[d, "\zeta^\Int"] \\
    f\Phi^\C_2 \mu^\C(x, \Phi^\C_2 \mu^\C(y, z)) \ar[r, "fa^\C"'] & f\Phi^\C_2\mu^\C(\Phi^\C_2 \mu^\C(x, y), z)
\end{tikzcd}
\end{equation*}

\begin{equation*}
\adjustbox{scale=0.92,center}{%
\begin{tikzcd}[column sep=4.5em]
    \Phi^\D_2\mu^\D(\Phi^\D_0(1^\D), f(x)) \ar[r, "{\Phi^\D_2\mu^\D(\nu^\Int, \id)}"] \ar[d, "\ell^\D"'] & \Phi^\D_2 \mu^\D(f\Phi^\C_0(1^\C), f(x)) \ar[d, "\zeta^\Int"] \\
    f(x) & f\Phi^\C_2\mu^\C(\Phi^\C_0(1^\C), x) \ar[l, "f \ell^\C"]
\end{tikzcd} \hspace{1.5em}
\begin{tikzcd}[column sep=4.5em]
    \Phi^\D\mu^2(f(x), \Phi^\D_0(1^\D)) \ar[r, "{\Phi^\D_2\mu^\D(\id, \nu^\Int)}"] \ar[d, "r^\D"'] & \Phi^\D_2 \mu^\D (f(x), f\Phi^\C_0(1^\C)) \ar[d, "\zeta^\Int"] \\
    f(x) & f\Phi^\C_2 \mu^\C(x, \Phi^\C_0(1^\C)) \ar[l, "fr^\C"]
\end{tikzcd}}
\end{equation*}where the (internalized) associator $a^\C$ and unitors $\ell^\C, r^\C$ for $\Int(\C_\sbull),$ and likewise $a^\D, \ell^\D,$ and $r^\D$ for $\Int(\D_\sbull),$ are as defined in \Cref{thm:internalization}.

Unfolding the definitions of each morphism in the hexagonal axiom, verifying commutativity amounts to checking that the perimeter of 
\begin{equation*}
\adjustbox{scale=0.85,center}{%
\begin{tikzcd}[column sep = 2.5em, row sep = 2.5em]
    \Phi^\D\mu^\D(f(x),\Phi^\D\mu^\D(f(y), f(z))) \ar[r, "\calR^\D"] \ar[d, "{\Phi^\D\mu^\D(\id, \Phi^\D\zeta)}"'] & \Phi^\D\mu^\D(f(x), \mu^\D(f(y), f(z))) \ar[r, "\Phi^\D\alpha^\D"] \ar[d, "{\Phi^\D\mu^\D(\id,\zeta)}"'] & \Phi^\D\mu^\D(\mu^\D(f(x), f(y)), f(z)) \ar[r, "(\calL^\D)^{-1}"] \ar[d, "{\Phi^\D\mu^\D(\zeta,\id)}"] & \Phi^\D\mu^\D(\Phi^\D\mu^\D(f(x), f(y)), f(z)) \ar[d, "{\Phi^\D\mu^\D(\Phi^\D\zeta, \id)}"] \\
    \Phi^\D\mu^\D(f(x), \Phi^\D f\mu^\D(f(y), f(z))) \ar[r, "\calR^\D"] \ar[d, "{\Phi^\D\mu^\D(\id, \xi)}"'] & \Phi^\D\mu^\D(f(x),f\mu^\C(y,z)) \ar[d, "\Phi^\D\zeta"'] & \Phi^\D\mu^\D(f\mu^\C(x,y),f(z)) \ar[r, "(\calL^\D)^{-1}"] \ar[d, "\Phi^\D\zeta"] & \Phi^\D\mu^\D(\Phi^\D f\mu^\C(x,y), f(z)) \ar[d, "{\Phi^\D\mu^\D(\xi, \id)}"] \\
    \Phi^\D\mu^\D(f(x),f\Phi^\C\mu^\C(y, z)) \ar[d, "\Phi^\D\zeta"'] & \Phi^\D f\mu^\C(x, \mu^\C(y,z)) \ar[r, "\Phi^\D f\alpha^\C"'] \ar[dd, "\xi"'] & \Phi^\D f\mu^\C(\mu^\C(x,y), z) \ar[dd, "\xi"] & \Phi^\D\mu^\D(f \Phi^\C\mu^\C(x,y), f(z)) \ar[d, "\Phi^\D\zeta"] \\
    \Phi^\D f\mu^\C(x,\Phi^\C\mu^\C(y, z)) \ar[d, "\xi"'] & & & \Phi^\D f\mu^\C(\Phi^\C\mu^\C(x, y), z) \ar[d, "\xi"] \\
    f\Phi^\C\mu^\C(x, \Phi^\C\mu^\C(y, z)) \ar[r, "f\calR^\C"'] & f\Phi^\C\mu^\C(x, \mu^\C(y, z)) \ar[r, "f\Phi^\C\alpha^\C"'] & f\Phi^\C\mu^\C(\mu^\C(x,y),z) \ar[r, "f(\calL^\C)^{-1}"'] & f\Phi^\C\mu^\C(\Phi^\C\mu^\C(x,y),z)
\end{tikzcd}
}
\end{equation*}
\noindent commutes, which follows from the commutativity of each inner sub-diagram:
\begin{itemize}[label=\textendash, noitemsep, topsep=2pt]
    \item the top left square commutes by naturality of $\calR^\D$;
    \item the bottom left heptagon is the coherence axiom \eqref{diag:lax-coh-right} of $\xi$ with the right absorbers;
    \item the top middle hexagon is the image of the associativity axiom of Definition \ref{def:gr-mon-functor} for the graded monoidal functor $f_\sbull$ under $\Phi^\D$;
    \item the bottom middle square commutes by naturality of $\xi$;
    \item the top right square commutes by naturality of $(\calL^{-1})^\D$;
    \item and the bottom right heptagon is the coherence axiom \eqref{diag:lax-coh-left} of $\xi$ with the left absorbers (with the absorbers inverted).
\end{itemize} Similarly, the compatibility of $(\zeta^\Int, \nu^\Int)$ with the left unitors is the outer perimeter of
\begin{equation*}
\adjustbox{scale=0.92,center}{%
\begin{tikzcd}
    \Phi^\D\mu^\D(\Phi^\D(1^\D), f(x)) \ar[rrr, "{\Phi^\D\mu^\D(\Phi^\D\nu,\id)}"] \ar[d, 
    "\calL^\D"'] &[-3em]&&[-3em] \Phi^\D\mu^\D(\Phi^\D f(1^\C),f(x)) \ar[rr, "{\Phi^\D\mu^\D(\xi,\id)}"] \ar[d, "\calL^\D"'] && \Phi^\D\mu^\D(f\Phi^\C(1^\C), f(x)) \ar[dd, "\Phi^\D\zeta"] \\
    \Phi^\D\mu^\D(1^\D,f(x)) \ar[rrr, "{\Phi^\D\mu^\D(\nu, \id)}"'] \ar[d, "\Phi^\D\lambda^\D"'] &&& \Phi^\D\mu^\D(f(1^\C), f(x)) \ar[d, "\Phi^\D\zeta"'] \\
    \Phi^\D f(x) \ar[d, "\iota^\D"'] &&& \ar[lll, "\Phi^\D f\lambda^\C"] \Phi^\D f\mu^\C(1^\C,x) \ar[dd, "\xi"'] && \Phi^\D f\mu^\C(\Phi^\C(1^\C), x) \ar[dd, "\xi"] \\
    f(x) & \Phi^\D f\Phi^\C(x) \ar[ul, "\Phi^\D f\iota^\C" {xshift=0.5em}] \ar[d, "\iota^\D"'] & |[xshift=5em]| \Phi^\D f\Phi^\C\mu^\C(1^\C,x) \ar[l, "\Phi^\D f\Phi^\C\lambda^\C"] \ar[ur, "\Phi^\D f\iota^\C"' {xshift=-0.5em}] \ar[dr, "\iota^\D"] \\
    & f\Phi^\C(x) \ar[ul, "f\iota^\C"] && f\Phi^\C\mu^\C(1^\C,x) \ar[ll, "f\Phi^\C\lambda^\C"] && f\Phi^\D\mu^\C(\Phi^\C(1^\C),x) \ar[ll, "f\calL^\C"]
\end{tikzcd}
}
\end{equation*}
wherein each interior sub-diagram commutes:
\begin{itemize}[label=\textendash, noitemsep, topsep=2pt]
    \item the top left square commutes by naturality of $\calL^\D$;
    \item the square below it is the image of left unitality axiom of Definition \ref{def:gr-mon-functor} for $f_\sbull$ under $\Phi^\D$;
    \item the trapezoid beneath that is the image of the naturality square of $\iota^\C$ under $\Phi^\D f$;
    \item the parallelogram at the bottom left commutes by the naturality of $\iota^\C$;
    \item the trapezoid to its left commutes by the naturality of $\iota^\D$;
    \item the triangle to the right of this is the coherence axiom \eqref{diag:lax-coh-triv} of $\xi$ with $\iota$;
    \item and the heptagon on the right is the coherence axiom \eqref{diag:lax-coh-left} of $\xi$ with the left absorbers.
\end{itemize} Commutativity of the remaining diagram, witnessing the compatibility of $(\zeta^\Int, \nu^\Int)$ with the right unitors, follows from a symmetric diagram chase. Hence we conclude that $\Int(f_\sbull) \colon \Int(\C_\sbull) \to \Int(\D_\sbull)$ is monoidal.
\end{proof}

Proposition \ref{prop:int_of_functor_is_monoidal} suggests a functorial interpretation of $\Int$, which we make formal shortly. First, we must detour briefly to describe the remaining 2-categorical structure on internalized categories.

Each internalized category $(\C_\sbull, \Phi^\C)$ is endowed with an internalized functor whose underlying graded monoidal functor is the identity, which we also denote simply by $\id_{\C_\sbull}.$ The structure maps $\xi$ for $\id_{\C_\sbull}$ are simply the identity componentwise, and the axioms are trivial. Furthermore, given a pair of internalized functors
\begin{equation*}
    (\C_\sbull, \Phi^\C) \xrightarrow{f_\sbull} (\D_\sbull, \Phi^\D) \xrightarrow{g_\sbull} (\E_\sbull, \Phi^\E)
\end{equation*}there is an internalized functor $(gf)_\sbull \colon (\C_\sbull, \Phi^\C) \to (\E_\sbull, \Phi^\E)$ whose underlying graded monoidal functor is the graded monoidal composite of $f_\sbull$ with $g_\sbull$. As one might guess, the structure maps are given by the composites
\begin{equation*}
    \Phi_m^\E g_mf_m \Rightarrow g_1\Phi^\D_m f_m \Rightarrow g_1f_1\Phi^\C_m,
\end{equation*} which again can be shown to satisfy the coherence conditions after some amount of tedious diagram chasing.

\begin{defn}
With the notions of identity and composition sketched above, we obtain a category $\IntMonCat$ of internalized graded monoidal categories and internalized functors between them.
\end{defn}

Of course, we are now compelled to define a suitable notion of natural transformation between internalized functors, to fulfill the prophecy that $\IntMonCat$ forms a 2-category.

\begin{defn} \label{def:int_trans}
Let $f_\sbull, g_\sbull \colon (\C_\sbull, \Phi^\C) \to (\D_\sbull, \Phi^\D)$ be internalized functors. An \emph{internalized (graded monoidal natural) transformation $\beta_\sbull \colon f_\sbull \Rightarrow g_\sbull$} is a graded monoidal transformation $f_\sbull \Rightarrow g_\sbull$ such that, for each $m \in \NN$, the following square commutes:
\begin{equation*}
\begin{tikzcd}
    \Phi^\D_mf_m \ar[r, Rightarrow, "\Phi_m^\D\beta_m"] \ar[d, Rightarrow, "\xi^m"'] & \Phi_m^\D g_m \ar[d, Rightarrow, "\xi^m"]\\
    f_1\Phi^\C_1 \ar[r, Rightarrow, "\beta_1\Phi^\C_1"'] & g_1\Phi_1^\C \period
\end{tikzcd}
\end{equation*}
\end{defn}

Just as internalized functors internalized to monoidal functors, internalized transformations likewise internalize to monoidal transformations:

\begin{prop} \label{prop:int_of_trans_is_monoidal}
Let $\beta_\sbull \colon f_\sbull \Rightarrow g_\sbull$ be an internalized transformation between internalized functors $f_\sbull, g_\sbull \colon \C_\sbull \to \D_\sbull$. Then $\Int(\beta) := \beta_1 \colon \Int(f) \Rightarrow \Int(g)$ defines a monoidal transformation between the monoidal functors.
\end{prop}
\begin{proof}
To show that $\beta_1$ is monoidal, we must verify that the diagrams
\begin{equation*}
\begin{tikzcd}
    & 1^\D \ar[dl, "\nu^\Int"'] \ar[dr, "\nu^\Int"]\\
    f_1\Phi^\C_0(1^\C) \ar[rr, "\beta_1"'] && g_1\Phi^\C_0(1^\C)
\end{tikzcd}
\end{equation*}and
\begin{equation*}
\begin{tikzcd}[column sep=4.5em]
    \Phi^\D_2\mu^\D(f_1(x), f_1(y)) \ar[r, "{\Phi^\D_2\mu^\D(\beta_1, \beta_1)}"] \ar[d, "\zeta^\Int"'] & \Phi_2^\D\mu^\D(g_1(x), g_1(y)) \ar[d, "\zeta^\Int"'] \\
    f_1\Phi^\C_2\mu^\C(x, y) \ar[r, "\beta_1"'] & g_1\Phi^\C_2\mu^\C(x, y)
\end{tikzcd}
\end{equation*}
commute (for all $x, y \in \C_1$). Unfolding the definitions from Proposition \ref{prop:int_of_functor_is_monoidal}, we can rewrite the first diagram as the perimeter of
\begin{equation}
\begin{tikzcd}
    & 1^\D \ar[dl, "\Phi^\D_0 \nu"'] \ar[dr, "\Phi^\D_0 \nu"]\\
    \Phi^\D_0f_0(1^\C) \ar[d, "\xi^0"'] \ar[rr, "\Phi^\D_0 \beta_0"'] && \Phi^\D_0g_0(1^\C) \ar[d, "\xi^0"]\\
    f_1\Phi^\C_0(1^\C) \ar[rr, "\beta_1"'] && g_1\Phi^\C_0(1^\C)
\end{tikzcd} \label{diag:int_trans_nu}
\end{equation}and the second as the perimeter of 
\begin{equation}
\begin{tikzcd}[column sep = 4.5em]
    \Phi^\D_2\mu^\D(f_1(x), f_1(y)) \ar[r, "{\Phi^\D_2\mu^\D(\beta_1, \beta_1)}"] \ar[d, "\Phi^\D_2 \zeta^{1,1}"'] & \Phi_2^\D\mu^\D(g_1(x), g_1(y)) \ar[d, "\Phi^\D_2 \zeta^{1,1}"] \\
    \Phi^\D_2 f_2 \mu^\C(x, y) \ar[r, "\Phi^\D_2 \beta_2"'] \ar[d, "\zeta^2"'] & \Phi^\D_2 g_2 \mu^\C(x, y)\ar[d, "\zeta^2"] \\
    f_1\Phi^\C_2\mu^\C(x, y) \ar[r, "\beta_1"'] & g_1\Phi^\C_2\mu^\C(x,y)
\end{tikzcd} \label{diag:int_trans_zeta}
\end{equation}Both diagrams commute because each of their interior subdiagrams commute:
\begin{enumerate}
    \item In diagram \eqref{diag:int_trans_nu} the top triangle is the image under $\Phi^\D_0$ of the coherence of $\beta_\sbull$ with the structure maps $\nu$, and the bottom square is the coherence of $\beta_\sbull$ with $\xi^0$.

    \item In diagram \eqref{diag:int_trans_zeta} the top square is the image under $\Phi^\D_2$ of the coherence of $\beta_\sbull$ with $\zeta^{1,1}$, and the bottom square is the coherence of $\beta_\sbull$ with $\xi^2$.
\end{enumerate}
It follows that $\Int(\beta)$ is a monoidal transformation $\Int(f) \Rightarrow \Int(g)$.
\end{proof}

Having defined $\Int(-)$ on internalized categories, functors, and natural transformations, we are finally ready to state the following:

\begin{prop} \label{prop:int-functorial}
Internalization defines a 2-functor $$\Int \colon \IntMonCat \to \MonCat,$$ which is the restriction of $\Forget_\Gr$ along the forgetful functors for monoidal categories and internalized categories. That is, the square of 2-categories
\begin{equation*}
\begin{tikzcd}
    \IntMonCat \ar[r, "\Int"] \ar[d, "\Forget_\Int"'] & \MonCat \ar[d, "\Forget_\Mon"] \\
    \GrMonCat \ar[r, "\Forget_\Gr"'] & \Cat
\end{tikzcd}
\end{equation*}commutes.
\end{prop}
\begin{proof}
The action of $\Int$ on objects is clear, while the actions on 1-morphisms and 2-morphisms is defined in Propositions \ref{prop:int_of_functor_is_monoidal} and \ref{prop:int_of_trans_is_monoidal} (where we already verified that these constructions satisfy the required conditions). Functoriality is automatic, since the composition and identities in $\IntMonCat$ restrict to classical functor composition and identities in degree 1.

Commutativity of the diagram is likewise straightforward: the underlying category (resp. functor, natural transformation) of an internalized category $\C_\sbull$ (resp. internalized functor $f_\sbull$, internalized transformation $\beta_\sbull$) is simply the degree 1 component $\C$ (resp. $f_1$, $\beta_1$). But the composite $\Forget_\Gr \circ \Forget_\Int$ likewise isolates the degree 1 data of an internalized category, by first forgetting the internalization and then concentrating in degree 1, so the desired square commutes.
\end{proof}

Next, we make formal the property that internalizing the free internalized category of a monoidal category recovers the original monoidal structure: 

\begin{thm} \label{thm:int-section}
The functor $\Int$ has a section $\Free \colon \MonCat \to \IntMonCat,$ which is the restriction of $\Free \colon \Cat \to \GrMonCat$ along the forgetful functors for monoidal categories and internalized categories. That is, the diagram of 2-categories
\begin{equation}
\begin{tikzcd}[column sep=4em, row sep=3.5em]
    \MonCat \ar[r, "\Free"] \ar[d, "\Forget_\Mon"'] & \IntMonCat \ar[d, "\Forget_\Int"] \\
    \Cat \ar[r, "\Free_\Gr"']  & \GrMonCat
\end{tikzcd} \label{diag:2-cat-square}
\end{equation}
commutes.
\end{thm}
\begin{proof}
For the diagram \eqref{diag:2-cat-square} to commute, for a monoidal category $(\C, \otimes, I)$ we are forced to take $\Free_\mathrm{Gr}(\C)$ to be the underlying graded monoidal category of $\Free(\C)$. We equip this graded monoidal category with the internalization $\Phi^\Free_\sbull,$ where $\Phi^\Free_0 \colon \mathbb{1} \to \C$ identifies the unit $I \in \C$, $\Phi^\Free_1$ is the identity functor on $\C$, and
\begin{align*}
    \Phi^\Free_m &\colon \C^m \to \C\\
    \Phi^\Free_m(x_1, \dots, x_m) &:= \bigotimes_{i=1}^m x_i
\end{align*}for all $m > 1$ (fixing a bias for $\otimes$, as usual). The absorbers of $\Phi^\Free_\sbull$ are induced by the monoidal structure of $\otimes$: if $m, n > 0$ then $\calL_{m, n}$ is the (formal) composite \cite{maclane} of the associators of $\otimes$ witnessing the isomorphism
\begin{align*}
    x_1 \otimes \cdots \otimes x_m \otimes \left(\bigotimes_{j = m+1}^{m+n} x_j\right) \cong \bigotimes_{i = 1}^{m + n} x_i
\end{align*} and $\calL_{0, n}$ is supplied via the left unitor
\begin{equation*}
    I \otimes \left(\bigotimes_{i = n}^n x_i\right) \cong \bigotimes_{i=1}^n x_i.
\end{equation*} $\calR$ is defined similarly, and by Mac Lane's Coherence Theorem these absorbers satisfy the axioms of \Cref{def:internalization}, making $(\Free(\C)_\sbull, \Phi^\Free_\sbull)$ an internalized category.

Given a monoidal functor $f \colon (\C, \otimes_\C, I_\C) \to (\D, \otimes_\D, I_\D),$ our definition of the underlying graded monoidal functor $\Free(f)_\sbull \colon \Free(\C)_\sbull \to \Free(\D)_\sbull$ is once again forced by \eqref{diag:2-cat-square},
\begin{equation*}
    \Free(f)_\sbull := \Free_\Gr(f)_\sbull.
\end{equation*} To extend this construction to a functor on (free) internalized categories, we must specify natural structure maps
\begin{equation*}
    \xi^m \colon \Phi^{\Free(\D)_\sbull}_m \Free(f)_m \Rightarrow f \Phi^{\Free(\C)_\sbull}_m,
\end{equation*}which, by unwinding the definitions, consist merely of maps
\begin{equation*}
    \bigotimes_{i=1}^m f(x_i) \to f\left(\bigotimes_{i=1}^m x_i\right)
\end{equation*} for each $x_1, \dots, x_m \in \C$ (where the empty tensor product is the unit). These are provided by the structure maps
\begin{align*}
    \zeta_{x, y} &\colon f(x) \otimes_\D f(y) \to f(x \otimes_\C y)\\
    \nu &\colon I_\D \to f(I_\C)
\end{align*} (where the maps $\zeta$ must be iterated consistently with the fixed choice of parenthesization for $\otimes$). The necessary diagrams commute by appealing to the Coherence Theorem (for lax monoidal functors -- see \cite[Section XI.2]{maclane}), using the axioms of $\zeta$ and $\nu$ as the structure maps of a monoidal functor to push the associators and unitors through them.

Finally, for a monoidal transformation $\beta \colon f \Rightarrow g$ between parallel monoidal functors, we again must take the graded monoidal transformation underlying $\Free(\beta)$ to be $\Free_\Gr(\beta)$. The condition of Definition \ref{def:int_trans} can be rewritten as
\begin{equation*}
\begin{tikzcd}
    \bigotimes_{i=1}^m f(x_i) \ar[r, "\otimes_i \beta"] \ar[d, "\xi^m"'] & \bigotimes_{i=1}^m g(x_i) \ar[d, "\xi^m"]\\
    f\left(\bigotimes_{i=1}^m x_i\right) \ar[r, "\beta"'] & g\left(\bigotimes_{i=1}^m x_i\right)
\end{tikzcd}
\end{equation*}commuting for each $x_1, \dots, x_m \in \C$. Since $\beta$ is a monoidal transformation and $\xi^m$ is defined using the structure maps $\zeta$ and $\nu$, each such square indeed commutes, completing the construction of $\Free \colon \MonCat \to \IntMonCat$. 2-functoriality is immediate, since the categorical structures on $\MonCat$ and $\IntMonCat$ are compatible and $\Free(-)$ is defined componentwise in every dimension.

Throughout our construction of $\Free(-)$ above, we took care to guarantee that the square \eqref{diag:2-cat-square} commutes, so it remains to be shown that $\Int \circ \Free = \id_\MonCat.$ Let $(\C, \otimes, I)$ be a monoidal category. On underlying categories we have
\begin{equation*}
\Int(\Free(\C)) = \Free(\C)_1 = \C^1 = \C,
\end{equation*} and given objects $x, y \in \C$ the internalized monoidal product of $x$ and $y$ is
\begin{equation*}
\Phi_2\mu_{1,1}(x, y) = \Phi_2(x, y) = x \otimes y.
\end{equation*} The same argument holds for morphisms in $\C,$ and the associator and unitors for the internalized product agree with those of $\otimes$ by the Coherence Theorem for $(\C, \otimes, I),$ since the the free monoidal product is strict, and the structure maps for the internalized product are formal composites of the structure maps of $\otimes.$ Thus, $\Free$ is indeed a right inverse for $\Int,$ completing the proof.
\end{proof}

By virtually the same proof as for graded monoidal categories, we obtain an adjunction $\Free \dashv \Int$ embedding monoidal categories in internalized categories.

\begin{thm} \label{thm:free-forget-adj-mon-int}
The functor $\Free$ is left adjoint to $\Int$, and embeds $\MonCat$ as a coreflective subcategory of $\IntMonCat.$
\end{thm}
\begin{proof}
The definitions of the unit $\eta \colon \id_\MonCat\Rightarrow \Int \circ \Free$ and $\epsilon \colon \Free \circ \Int \Rightarrow \id_\IntMonCat$ are identical to the constructions of Proposition \ref{prop:free-forget-adj-cat-grmoncat}, and since $\Int \circ \Free = \id_\MonCat$ by Theorem \ref{thm:int-section} the same argument shows that this lifts to an adjunction $\Free \dashv \Int$. As before, since $\Int$ is a retraction of $\Free$ the latter is an embedding, and since $\eta$ is the identity componentwise it follows that $\Free$ embeds $\MonCat$ as a coreflective subcategory.
\end{proof}

\begin{cor}
Internalization preserves all limits in $\IntMonCat$ (as a 1-category).
\end{cor}
\begin{proof}
Immediate from \Cref{thm:free-forget-adj-mon-int}.
\end{proof}

\section{Internal Smash Products of Spectra}
\label{section:spectra}

We conclude this paper by outlining two applications of internalized categories to reconstruct a pair of monoidal products which occur naturally in stable homotopy theory. 

\subsection{Day Convolution}
\label{subsection:day_convolution}

The smash product of symmetric or orthogonal spectra (or, more generally, of diagram spectra over a given monoid object in the category of $\D$-spaces, where $\D$ is some topological category -- see \cite{MMSS}) can be formulated as a special case of Day convolution. In this section, we describe how Day convolution fits into the framework of an internalized graded monoidal category, and use this to prove that the Day convolution monoidal product is closed and symmetric.

We first recall the general procedure of Day convolution. We fix a bicomplete, closed, symmetric monoidal category $(\V, \otimes_\V, I_\V)$, and a (small) $\V$-enriched monoidal category $(\C, \otimes_\C, I_\C)$. We additionally denote by $\Fun(\C, \V)$ the $\V$-enriched functor category, with hom-object $\Fun(\C, \V)(F, G)$ defined by the $\V$-enriched end
\begin{equation*}
    \Fun(\C, \V)(F, G) := \int_{c \in \C} \underline{V}(Fc, Gc)
\end{equation*} for given $\V$-functors $F, G \colon \C \to \V$. Readers unfamiliar with this notation or who wish to review (co)ends should consult Section 1.2 of \cite{riehl} or Section 1.1.3 of \cite{coend}. Lastly, given $\V$-functors $F: \C \to \A$ and $K: \C \to \B$ (where $\A$ and $\B$ are $\V$-categories), we write $\Lan_K(F)$ for the left Kan extension $\B \to \A$ of $F$ along $K$, as defined in Section 4.1 of \cite{kelly}.

\begin{defn} \label{def:day-conv-prod}
The \emph{Day convolution product} on $\Fun(\C, \V)$ is the bifunctor
\begin{equation*}
    \Daytimes \colon \Fun(\C, \V) \times \Fun(\C, \V) \to \Fun(\C, \V)
\end{equation*} defined by the left Kan extension
\begin{equation*}
    F \Daytimes G := \Lan_{\otimes_\C}(F \otimes_\V G)
\end{equation*} for $\V$-functors $F, G \in \Fun(\C, \V)$, where
\begin{equation*}
    (F \otimes_\V G)(c) := Fc \otimes_\V Gc.
\end{equation*} By way of the coend formula for left Kan extensions (see (4.25) in \cite{kelly}), we compute the image of $c \in \C$ under $F \Daytimes G$ as the coend in $\V$
\begin{equation*}
    (F \Daytimes G)(c) = \int^{c_1, c_2 \in \C} \C(c_1 \otimes_\C c_2, c) \otimes_\V Fc_1 \otimes_\V Gc_2.
\end{equation*}
\end{defn}

\begin{rem}
The product $\Daytimes$ is indeed monoidal, with unit object given by $\Y(I_\C) \colon \C \to \V$, where $\Y$ is the $\V$-enriched (contravariant) Yoneda embedding $\C^\op \to \Fun(\C, \V)$. One can prove this classically, but this will also follow immediately from \Cref{thm:internalization} after we've shown that $\Daytimes$ arises as the internalization of a graded monoidal product.
\end{rem}

\begin{defn}
We define a $\V$-enriched graded monoidal category $(\Fun_\sbull(\C, \V), \mu, 1)$ by
\begin{enumerate}
    \item $\Fun_m(\C, \V) := \Fun(\C^m, \V)$ for $m \in \NN$.
    \item $\mu_{m, n}(F, G) \colon \C^{m + n} \to \V$ is the $\V$-functor
    \begin{equation*}
        \mu_{m, n}(F, G)(c_1, \dots, c_{m + n}) := F(c_1, \dots, c_m) \otimes_V G(c_{m + 1}, \dots, c_{m + n})
    \end{equation*} for $m, n \in \NN$, where $F \colon \C^m \to \V$ and $G \colon \C^n \to \V$ are $\V$-functors.
    \item $1 \colon \C^0 \to \V$ is the functor which carries the terminal $\V$-category $\C^0$ to the object $I_\V \in \V$.
\end{enumerate} The natural isomorphisms for $\mu$ are defined termwise by the corresponding isomorphisms for $\otimes_\V$.
\end{defn}

We aim to show that $\Daytimes$ is an internalization of $\mu$.

\begin{defn}
Given $m \in \NN$, let $\Phi_m \colon \Fun(\C^m, \V) \to \Fun(\C, \V)$ be the $\V$-functor
\begin{equation*}
    \Phi_m(F) := \Lan_{\otimes_\C^m}(F),
\end{equation*} where $\otimes_\C^m \colon \C^m \to \C$ is the\footnote{Here one should be slightly careful in the case that $\otimes_\C$ is not strict. One is then forced to choose a particular arrangement of parentheses, that is, a bias for the $m$-fold monoidal product on $\C$. Each natural isomorphism in the internalization is then composed with some sequence of associators to ensure compatibility, but Mac Lane's Coherence Theorem guarantees that this will not cause problems.} $\V$-functor
\begin{equation*}
    (c_1, \dots, c_m) \mapsto \bigotimes_{i = 1}^m c_i.
\end{equation*} In the special case that $m = 0$, we adopt the convention that $\otimes_\C^0$ is the $\V$-functor $\C^0 \to \C$ which identifies the unit $I_\C$.
\end{defn}

\begin{rem} \label{rem:0_is_yoneda}
We note that since $\Fun_0(\C, \V) = \Fun(\C^0, \V)$ is naturally isomorphic to $\V$, $\Phi_0$ can be viewed as a $\V$-functor $\V \to \Fun(\C, \V)$. Given $v \in \V$ and $c \in \C$, we have
\begin{equation*}
    \Phi_0(v)(c) = (\Lan_{\otimes_\C^0}(v))(c) = \int^{\C^0} \C(I_\C, c) \otimes_\V v = \Y(I_\C)(c) \otimes_\V v
\end{equation*} by the coend formula for left Kan extensions.
\end{rem}

\begin{thm}
The functor $\Phi$ defines an internalization on $(\Fun_\sbull(\C, \V), \mu, 1)$.
\end{thm}
\begin{proof}
First, we observe that $\Phi_1(F)$ is the left Kan extension of $F$ along $\id_\C$, so there is a natural isomorphism $u \colon \Phi_1 \cong \id_{\Fun(\C, \V)}$. 

Next, we will construct the left absorber for $\Phi$ as the composite of a chain of natural isomorphisms; the construction of the right absorber is completely dual.

Fix $m, n \in \NN$, and let $F \colon \C^m \to \V$ and $G \colon \C^n \to \V$ be $\V$-functors. Additionally, let $c \in \C$. We aim to construct a (natural) isomorphism
\begin{equation*}
    \Lan_{\otimes_\C^{n + 1}}(\Lan_{\otimes^m_\C}(F), G)(c) \cong \Lan_{\otimes_\C^{m + n}}(F, G)(c).
\end{equation*} Unwinding the left-hand side using the coend formula, we obtain
\begin{equation*}
    \int^{d, c_1, \dots, c_n} \C\Big(d \otimes_\C \bigotimes_{i = 1}^n c_i, c\Big) \otimes_\V \left(\int^{d_1, \dots, d_m} \C\Big(\bigotimes_{j = 1}^m d_j, d\Big) \otimes_\V F(d_1, \dots, d_m)\right) \otimes_\V \hspace{0.2em}G(c_1, \dots, c_n),
\end{equation*} where both coends are indexed by $\C$. By the Fubini Theorem, this is naturally isomorphic to
\begin{equation*}
    \int^{d, d_1, \dots, d_m, c_1, \dots, c_n} \C\Big(d \otimes_\C \bigotimes_{i = 1}^n c_i, c\Big) \otimes_\V \hspace{0.2em} \C\Big(\bigotimes_{j = 1}^m d_j, d\Big) \otimes_\V\hspace{0.2em} F(d_1, \dots, d_m) \otimes_\V \hspace{0.2em}G(c_1, \dots, c_n),
\end{equation*} which is naturally isomorphic to
\begin{equation*}
    \int^{d_1, \dots, d_m, c_1, \dots, c_n} \C\Big(\bigotimes_{j = 1}^m d_j \otimes_\C \bigotimes_{i = 1}^n c_i, c\Big)  \otimes_\V \hspace{0.2em} F(d_1, \dots, d_m) \otimes_\V \hspace{0.2em}G(c_1, \dots, c_n)
\end{equation*} via the co-Yoneda Lemma. But using the definition of $\Phi$ and the coend formula for left Kan extensions, this is precisely
\begin{equation*}
    \Lan_{\otimes_\C^{m + n}}(F, G)(c),
\end{equation*} so the desired natural isomorphism is the composite of the isomorphisms above. Note that in the case that one of $m$ or $n$ is $0$ the same argument holds, with the corresponding $m$-fold or $n$-fold tensor products being replaced by the unit object $I_\C$.

For brevity we omit the verification of the coherence axioms for the absorbers.
\end{proof}

The following result is then immediate from \Cref{thm:internalization} and \Cref{def:day-conv-prod}.

\begin{prop}\label{prop:day_is_internalization}
The monoidal product on $\Fun(\C, \V)$ induced by $\Phi$ is the Day convolution product $\Daytimes$.
\end{prop}

\begin{cor}
The Day convolution product is a monoidal product, with unit $\Y(I_\C)$.
\end{cor}
\begin{proof}
It follows from \Cref{prop:day_is_internalization} and \Cref{thm:internalization} that $\Daytimes$ is a monoidal product, while the second claim is a consequence of \Cref{rem:0_is_yoneda}, since
\begin{equation*}
    \Phi_0(1)(c) = \Y(I_\C)(c) \otimes_\V I_\V \cong \Y(I_\C)(c)
\end{equation*} for each $c \in \C$.
\end{proof}

We have shown that Day convolution is a special case of internalization. In fact, we can further recover the properties of symmetry and closure that $\Daytimes$ enjoys by using the formal behavior of internalized categories.

\begin{lem}
The graded monoidal category $\Fun_\sbull(\C, \V)$ is symmetric, as is its internalization (see \Cref{def:symmetric_monoidal} and \Cref{def:symmetric_internalization}).
\end{lem}
\begin{proof}
This follows from the symmetry of $\otimes_\V$, and a diagram chase on its braiding and the absorbers of $\Phi$.
\end{proof}

\begin{lem} \label{lem:Day_closed}
The graded monoidal category $\Fun_\sbull(\C, \V)$ is closed (see \Cref{def:closed_monoidal}).
\end{lem}
\begin{proof}
We define the internal hom of $\mu$ by the $\V$-functors
\begin{align*}
    \scF_{m, n} &\colon \Fun(\C^n, \V)^\op \times \Fun(\C^{m + n}, \V) \to \Fun(\C^m, \V)
\end{align*} by
\begin{equation*}
    \scF_{m, n}(F, G)(c_1, \dots, c_m) := \int_{d_1, \dots, d_n \in \C} \underline{V}(F(d_1, \dots, d_n), G(c_1, \dots, c_m, d_1, \dots, d_n)),
\end{equation*} where $c_1, \dots, c_m \in \C$, $F \colon \C^n \to \V$, and $G \colon \C^{m + n} \to \V$. To verify that this does indeed define an internal hom of $\mu$, fix $m, n \in \NN$, and suppose given $\V$-functors $F \colon \C^m \to \V$, $G \colon \C^n \to \V$, and $H \colon \C^{m + n} \to \V$. Then, unwinding the definitions, we have
\begin{equation*}
    \Fun_{m + n}(\C, \V)(\mu(F, G), H) = \int_{c_1, \dots, c_{m + n} \in \C} \hspace{-4em}\underline{V}(F(c_1, \dots, c_m) \otimes_\V G(c_{m+1}, \dots, c_{m+n}), H(c_1, \dots, c_{m+n})),
\end{equation*} and the right-hand side is isomorphic to
\begin{equation*}
    \int_{c_1, \dots, c_{m + n} \in \C} \underline{V}(F(c_1, \dots, c_m), \underline{V}(G(c_{m + 1}, \dots, c_{m + n}), H(c_1, \dots, c_{m + n})))
\end{equation*} using the closed structure of $\V$. Now, using the universal property of ends, this is isomorphic to
\begin{equation*}
    \int_{c_1, \dots, c_m \in \C} \underline{V}(F(c_1, \dots, c_m), \int_{c_{m + 1, \dots, c_{m + n}}} \underline{V}(G(c_{m+1}, \dots, c_{m+n}), H(c_1, \dots, c_{m+n}))),
\end{equation*} which is precisely
\begin{equation*}
    \Fun_m(\C, \V)(F, \scF(G, H)),
\end{equation*} and hence $\mu$ is closed with internal hom $\scF$.
\end{proof}

\begin{lem}\label{lem:Day_internalization_closed}
The internalization $\Phi$ of $\mu$ has a right adjoint, and is closed (see \Cref{def:closed_internalization}).
\end{lem}
\begin{proof}
Since $\mu$ is closed by \Cref{lem:Day_closed}, the second claim follows immediately from the first. Given $m \in \NN$ and a $\V$-functor $G \colon \C^m \to \V$, we define a $\V$-functor $\Psi_m \colon \Fun(\C, \V) \to \Fun(\C^m, \V)$ by
\begin{equation*}
    \Psi_m(G)(c_1, \dots, c_m) := G\Big(\bigotimes_{i = 1}^m c_i\Big)
\end{equation*} for each $c_1, \dots, c_m \in \C$. Then since $\Phi_m$ is given by left Kan extension along $\otimes_\C^m$, and $\Psi_m$ precomposes with $\otimes_\C^m$, it is automatic that $\Phi_m$ is left adjoint to $\Psi_m$.
\end{proof}

\begin{prop}
The Day convolution product is symmetric and closed.
\end{prop}
\begin{proof}
Symmetry follows from \Cref{cor:symmetric_internalization} and \Cref{lem:Day_closed}, and closure follows from \Cref{cor:closed_internalization}, \Cref{lem:Day_closed}, and \Cref{lem:Day_internalization_closed}.
\end{proof}

To specialize to the case of spectra (symmetric, orthogonal, etc.) we refer to \cite{MMSS}. Fix a symmetric monoidal category $\D$ enriched over based compactly generated topological spaces $\Top_*$, and recall that a \emph{$\D$-space} is a continuous functor $\D \to \Top_*$. Taking $\V := \Top_*$, with symmetric monoidal product given by the usual smash product $\wedge$, and $\C := \D$, we see that $\Fun(\C, \V)$ is precisely the category $\D(\Top_*)$ of $\D$-spaces. Day convolution then provides a means of turning the external product induced by $\wedge$ into the \emph{internal smash product of $\D$-spaces} -- see Definition 21.4 of \cite{MMSS}.

A \emph{$\D$-spectrum} (Definition 1.9 in \cite{MMSS}) over a monoid object $R$ in $\D(\Top_*)$ is a $\D$-space equipped with a right action by $R$. These form a category $\categ{Spec}(\D)_R$ which, by Theorem 2.2 of \cite{MMSS}, is equivalent to the category of $\D_R$-spaces, where $\D_R$ is a $\Top_*$-enriched category defined from $\D$ and $R$ (see Construction 2.1 of \cite{MMSS} for a complete definition). Moreover, when $R$ is a commutative monoid $\D_R$ is naturally equipped with a symmetric monoidal product, in which case Day convolution induces a symmetric monoidal product on the category of $\D$-spectra.

Finally, we note that both symmetric spectra and orthogonal spectra can be realized as $\D$-spectra over commutative monoids. In the case of the former (see Example 4.2 of \cite{MMSS}), we have $\D := \Sigma$, the category of (unbased) finite sets $[n] := \{0, \dots, n\}$ for each $n \in \NN$ and their permutations, and $R := S_\Sigma$, the continuous functor $\Sigma \to \Top_*$ which sends $[n]$ to the $n$-sphere $S^n$. For the latter (Example 4.4 of \cite{MMSS}) we instead take $\D$ to be the category of (unbased) finite-dimensional real inner product spaces and isometric linear isomorphisms between them, and let $R$ be the functor $\D \to \Top_*$ sending $V$ to the one-point compactification $S^V$. In both cases, the resulting internal smash product induced by Day convolution is precisely the usual smash product of the corresponding spectra.

\subsection{The Smash Product of \textit{S}-Modules}
\label{subsection:s-modules}

Finally, we turn our attention to $S$-modules, defined in \cite{EKMM}. After reviewing the relevant definitions, we will show how the smash product of $S$-modules can be realized as an internalization of the external smash product of coordinate-free spectra.

Recall that, given a real inner product $U$ isomorphic to $\RR^\infty$, the direct sum of countably-many copies of $\RR$, there exists a category $\scS U$ of \emph{coordinate-free spectra over the universe $U$}. Throughout the remainder of this section, we will refer to these as ``spectra (over $U$)."

The category $\scS U$ enjoys a number of convenient properties. By Proposition I.1.1 in \cite{EKMM} it is both complete and cocomplete. Moreover, given a pair of universes $U$ and $U'$, there is an \emph{external smash product}
\begin{equation*}
    \wedge \colon \scS U \times \scS U' \to \scS( U \oplus U')
\end{equation*}
of spectra over $U$ and $U'$, which is associative and symmetric (up to natural isomorphism). Lastly, there is a functor $\Sigma^\infty \colon \Top_* \to \scS U$, called the \emph{suspension spectrum functor}, which identifies a based space with a spectrum over $U$, and we write $S$ to denote the suspension spectrum $\Sigma^\infty S^0$.

For universes $U$ and $U'$, let $\scI(U, U')$ be the set of linear isometries $U \to U'$. Viewing $U$ and $U'$ as topological spaces, whose topologies are induced by their finite dimensional subspaces, we can equip $\scI(U, U')$ with the function space topology. Now, given any space $A$ and a continuous map $\alpha \colon A \to \scI(U, U')$, together with a spectrum $E$ over $U$, there exists a spectrum $A \ltimes E$ over $U'$, called the \emph{twisted half-smash product} of $A$ and $E$. Moreover, this construction is functorial in $E$, and satisfies the following properties (see Proposition I.2.2 of \cite{EKMM}).
\begin{enumerate}
    \item \label{item:ltimes_id} For each $E \in \scS U$ there is a canonical isomorphism of spectra
    \begin{equation*}
        \{\id_U\} \ltimes E \cong E.
    \end{equation*}
    \item \label{item:ltimes_assoc} For spaces $A, B$ with maps $\alpha \colon A \to \scI(U, U')$ and $\beta \colon B \to \scI(U', U'')$ and a spectrum $E \in \scS U$ there is a canonical isomorphism
    \begin{equation*}
        (B \times A) \ltimes E \cong B \ltimes (A \ltimes E).
    \end{equation*}
    \item \label{item:ltimes_distribute} Let $A_1$ and $A_2$ be spaces with maps $\alpha_i \colon A_i \to \scI(U_i, U_i')$. Then given spectra $E_1 \in \scS U_1$ and $E_2 \in \scS U_2$ there is a canonical isomorphism
    \begin{equation*}
        (A_1 \times A_2) \ltimes (E_1 \wedge E_2) \cong (A_1 \ltimes E_1) \wedge (A_2 \ltimes E_2).
    \end{equation*}
    \item \label{item:ltimes_homo} For each $E \in \scS U$ and each based space $X$ there is a canonical isomorphism
    \begin{equation*}
        A \ltimes (E \wedge X) \cong (A \ltimes E) \wedge X.
    \end{equation*}
\end{enumerate}

Now, writing $U^j$ for the direct sum of $j$ copies of a universe $U$, and taking $\scL(j) := \scI(U^j, U)$, there is a (symmetric) operad, called the \emph{linear isometries operad}, equipped with structure maps
\begin{equation*}
    \gamma \colon \scL(k) \times \scL(j_1) \times \cdots \times \scL(j_k) \to \scL(j_1 + \cdots + j_k),
\end{equation*} defined by
\begin{equation*}
    \gamma(g; f_1, \dots, f_k) := g (f_1, \dots, f_k).
\end{equation*} For each $j$ and each spectrum $E \in \scS U^j$, the identity map on $\scL(j)$ then induces a spectrum $\scL(j) \ltimes E \in \scS U$. (In the special case that $j = 0$, we in fact have that $\scL(0) \ltimes X \cong \Sigma^\infty X$ for any based topological space $X$.) Moreover, the endofunctor
\begin{equation*}
    \scL(1) \ltimes (-) \colon \scS U \to \scS U
\end{equation*} defines a monad, which we denote by $\LL$, with multiplication map $\LL \LL E \to \LL E$ induced by the structure map $\gamma \colon \scL(1) \times \scL(1) \to \scL(1)$, and unit $E \to \LL E$ induced by the inclusion $\{0\} = \scL(0) \to \scL(1)$ which identifies $\id_U$. 

We write $\scS[\LL]$ for the category of $\LL$-algebras, which we refer to as \emph{$\LL$-spectra}. This category is equipped with a bifunctor $\wedge_\scL \colon \scS[\LL] \times \scS[\LL] \to \scS[\LL]$, called the \emph{smash product of $\LL$-spectra}, defined as follows: let $M$ and $N$ be $\LL$-spectra, with actions $\xi_M \colon \LL M \to M$ and $\xi_N \colon \LL N \to N$; we take $M \wedge_\scL N $ to be the coequalizer of the diagram of $\LL$-algebras
\begin{equation*}
    \begin{tikzcd}[column sep = 4em]
        (\scL(2) \times \scL(1) \times \scL(1)) \ltimes (M \wedge N) \ar[r, shift right = 0.35em, "\gamma \ltimes \id"'] \ar[r, shift left = 0.35em, "\id \ltimes (\xi_M \wedge \xi_N)"] & \scL(2) \ltimes (M \wedge N).
    \end{tikzcd}
\end{equation*} 

\begin{notn}
    Following the notation of \cite{EKMM}, we denote the coequalizer of
    \begin{equation*}
    \begin{tikzcd}
        (\scL(j_1) \times \cdots \times \scL(j_n) \times \scL(k) \times \scL(1)^k) \ltimes E \ar[r, shift right = 0.35em, "(\id \times \gamma) \ltimes \id"'] \ar[r, shift left = 0.35em, "\id \ltimes \xi_E"] & (\scL(i_1) \times \cdots \times \scL(i_n) \times \scL(k)) \ltimes E
    \end{tikzcd}
    \end{equation*}by
    \begin{equation*}
        (\scL(j_1) \times \cdots \times \scL(j_n) \times \scL(k)) \ltimes_{\scL(1)^k} E.
    \end{equation*}
\end{notn}

As shown in Proposition I.5.2 and Theorem I.5.5 of \cite{EKMM} the smash product $\wedge_\scL$ is symmetric and associative (up to natural isomorphism), and for each $M \in \scS [\LL]$ there is a canonical map
\begin{equation*}
    S \wedge_\scL M \to M,
\end{equation*} described in Proposition I.8.3 of \cite{EKMM}.

An \emph{$S$-module} is an $\LL$-spectrum for which this map is an isomorphism, and we write $\scM_S$ for the full subcategory of $\scS[\LL]$ spanned by the $S$-modules. For $S$-modules $M$ and $N$, we write $M \wedge_S N := M \wedge_\scL N$, which is again an $S$-module (see Proposition II.1.2 of \cite{EKMM}), so $\wedge_S$ defines a symmetric monoidal product on $\scM_S$.

We aim to show that $\wedge_S$ is realized as the internalization of the external product of spectra $\wedge$. First, we must formulate a suitable graded monoidal category to impose an internalization upon. This is somewhat tricky, primarily because of the non-unitality of $\wedge_\scL$. We essentially need to define a kind of intermediate, ``non-unital" internalized category (i.e. we forget about the degree $0$ component, and only require an associative product) which acts like $\scS[\LL]$ above, and then restrict to those objects which behave like $S$-modules to recover unitality.

\begin{defn}
    Fix $m \geq 1$. Then there is a monad $\LL^m$ on $\scS U^m$, defined by
    \begin{equation*}
        \LL^m(E) := \scL(1)^m \ltimes E,
    \end{equation*} where the map $\scL(1)^m \to \scI(U^m, U^m)$ carries a tuple of isometries $(f_1, \dots, f_m)$ to the isometry $f_1 \times \cdots \times f_m$. The multiplication map $\mu_m \colon \LL^m \LL^m E \to \LL^m E$ is induced by $m$ copies of the structure map $\gamma \colon \scL(1) \times \scL(1) \to \scL(1)$, while the unit $\eta_m :\colon E \to \LL^m E$ is induced by the inclusion of $\scL(0)$ into $\scL(1)^m$ which identifies $(\id_U, \dots, \id_U)$.
    
    The category $\scS[\LL^m]$ of \emph{$\LL^m$-spectra} is the category of algebras over $\LL^m$.
\end{defn}

Setting $\scS[\LL]_0 := \Top_*$, the categories $\scS[\LL]_m$ assemble to form a graded category $\scS[\LL]_\sbull$. The external smash product of spectra then defines an associative graded product on $\scS[\LL]_\sbull$.

\begin{defn}
    For $m, n \geq 1$, there is a bifunctor $\mu_{m, n} \colon \scS[\LL]_m \times \scS[\LL]_n \to \scS[\LL]_{m + n}$, which, given $(M, \xi_M) \in \scS[\LL]_m$ and $(N, \xi_N) \in \scS[\LL]_n$, equips $M \wedge N$ with the structure map $\xi \colon \LL^{m + n}(M \wedge N) \to M \wedge N$ defined by the composite
    \begin{equation*}
        (\scL(1)^m \times \scL(1)^n) \ltimes (M \wedge N) \cong (\scL(1)^m \ltimes M) \wedge (\scL(1)^n \ltimes N) \xrightarrow{\xi_M \wedge \xi_N} M \wedge N.
    \end{equation*} 

    If $n = 0$ then we take $\mu_{m, 0} \colon \scS[\LL]_m \times \scS[\LL]_0 \scS[\LL]_m$ to be the smash product of an $\LL^m$-spectrum $M$ with a (based) space $X$. That is,
    $$
        \mu_{m, 0}(M, X) := M \wedge X,
    $$
    which is an $\LL^m$-spectrum by the isomorphism of property \ref{item:ltimes_homo}. Similarly, if $m = 0$ then $\mu_{0, n}(X, N) := X \wedge N$, and if $m, n = 0$ then $\mu_{0, 0}$ is the ordinary smash product of based spaces.
\end{defn}

The category $\scS[\LL]_\sbull$ will lay the foundation for our internalization of the smash product of spectra. However, before we can begin assembling this internalization, we will need a preliminary result, which generalizes Lemma I.5.4 of \cite{EKMM}.

\begin{lem}[Hopkins] \label{lem:hopkins}
For integers $k \geq 1$, $j_1, \dots, j_k \geq 1$, the diagram

\adjustbox{scale=0.9,center}{%
    \begin{tikzcd}
        \scL(k) \times \scL(1)^k \times \scL(j_1) \times \cdots \times \scL(j_k) \ar[r, shift right = 0.35em, "\gamma \times \id"'] \ar[r, shift left = 0.35em, "\id \times \gamma^k"] & \scL(k) \times \scL(j_1) \times \cdots \times \scL(j_k) \ar[r, "\gamma"] & \scL(j_1 + \cdots + j_k)
    \end{tikzcd}
} is a split coequalizer of spaces, and consequently
\begin{equation*}
    \scL(j_1 + \cdots + j_k) \cong \scL(k) \times_{\scL(1)^k} \left(\scL(j_1) \times \cdots \times \scL(j_k)\right).
\end{equation*}
\end{lem}
\begin{proof} The argument is a direct generalization of the one presented in \cite{EKMM}. For each $i = 1, \dots, k$ we fix an isomorphism $s_i \colon U^{i_k} \cong U$. We define maps
\begin{align*}
    h \colon \scL(j_1 + \cdots + j_k) &\to \scL(k) \times \scL(j_1) \times \cdots \times \scL(j_k)\\
    f &\mapsto (f \Big(\bigoplus_{i = 1}^k s_i\Big)^{-1}; s_1, \dots, s_k)
\end{align*} and
\begin{align*}
    r \colon \scL(k) \times \scL(j_1) \times \cdots \times \scL(j_k) &\to \scL(k) \times \scL(1)^k \times \scL(j_1) \times \cdots \times \scL(j_k)\\
    (f; g_1, \dots, g_k) &\mapsto (f; g_1 s_1^{-1}, \dots, g_k s_k^{-1}; s_1, \dots, s_k),
\end{align*} and claim that these split the fork above.

First, given $f \in \scL(j_1 + \cdots + j_k)$, we observe that
\begin{equation*}
    \gamma h(f) = f  \Big(\bigoplus_{i = 1}^k s_i\Big)^{-1}(s_1, \dots, s_k) = f \bigoplus_{i = 1}^k \id_{U^{j_i}} = f,
\end{equation*} so $\gamma h$ is the identity. Additionally, given $(f; g_1, \dots, g_k) \in \scL(k) \times \scL(j_1) \times \cdots \times \scL(j_k)$, we have
\begin{align*}
    (\id \times \gamma^k)r(f; g_1, \dots, g_k) &= (\id \times \gamma^k)(f; g_1 s_1^{-1}, \dots, g_k s_k^{-1}; s_1, \dots, s_k)\\
    &= (f; g_1, \dots, g_k),
\end{align*} so $(\id \times \gamma^k)r$ is the identity as well. Finally, 
\begin{align*}
    (\gamma \times \id)r(f; g_1, \dots, g_k) &= (\gamma \times \id)(f; g_1 s_1^{-1}, \dots, g_k s_k^{-1}; s_1, \dots, s_k)\\
    &= (f(g_1s^{-1}, \dots, g_ks_k^{-1}); s_1, \dots, s_k),
\end{align*} and
\begin{align*}
    h\gamma (f; g_1, \dots, g_k) &= h (f(g_1, \dots, g_k))\\
    &= (f(g_1, \dots, g_k) \Big(\bigoplus_{i = 1}^k s_i\Big)^{-1}; s_1, \dots, s_k)\\
    &= (f(g_1s^{-1}, \dots, g_ks_k^{-1}); s_1, \dots, s_k),
\end{align*} so we conclude that $(\gamma \times \id)r = h\gamma$.
\end{proof}

We now begin to construct the desired internalization $\Phi_\sbull$. It may seem strange to do this before defining the graded monoidal category that we ultimately aim to internalize, but this category will in fact be defined using our ``incomplete" internalization of $\scS[\LL]_\sbull$.

\begin{defn}
    Given $M \in \scS[\LL]_m$, we define the underlying object of $\Phi_m(M, \xi_M)$ to be the coequalizer of the diagram in $\scS U$
    \begin{equation*}
        \begin{tikzcd}
            (\scL(m) \times \scL(1)^m) \ltimes M \ar[r, shift right = 0.35em, "\gamma \ltimes \id"'] \ar[r, shift left = 0.35em, "\id \ltimes \xi_M"] & \scL(m) \ltimes M.
        \end{tikzcd}
    \end{equation*} The structure map $\LL \Phi_m(M) \to \Phi_m(M)$ is induced by $\gamma \colon \scL(1) \times \scL(m) \to \scL(m)$. 
\end{defn}

We would hope that the collection of functors $\Phi_m$ would assemble to form an internalization $\Phi$ on the graded monoidal category $(\scS[\LL]_\sbull, \mu, S^0)$. Indeed, this is very nearly the case: in all higher (i.e. nonzero) degrees, these functors do in fact enjoy both left and right absorbers.

\begin{lem} \label{lem:higher_absorbers}
    Let $m, n > 0$, let $M \in \scS[\LL]_m$, and let $N \in \scS[\LL]_n$. Then there exist natural isomorphisms
    \begin{align*}
        \Phi_{n + 1}\mu(\Phi_m M, N) \cong \Phi_{m + n}\mu(M, N) \cong \Phi_{m + 1}\mu(M, \Phi_n N).
    \end{align*}
\end{lem}
\begin{proof}
    We will show explicitly only that the left isomorphism exists, as the right isomorphism follows from a symmetric argument. Since $N$ is an $\LL^n$-algebra there is a split coequalizer diagram in $\scS U$.
\begin{equation*}
    \begin{tikzcd}
        \LL^n \LL^n N \ar[r, shift right = 0.25em] \ar[r, shift left = 0.25em] & \LL^n N \ar[r] & N,
    \end{tikzcd}
\end{equation*} so $N \cong \scL(1)^n \ltimes_{\scL(1)^n} N$. Thus, we have
\begin{align*}
    \Phi_{n + 1}\mu(\Phi_m M, N) &= \scL(n + 1) \ltimes_{\scL(1)^{n + 1}} (\scL(m) \ltimes_{\scL(1)^m} M \wedge N)\\
    &\cong \scL(n + 1) \ltimes_{\scL(1)^{n + 1}} (\scL(m) \ltimes_{\scL(1)^m} M \wedge \scL(1)^n \ltimes_{\scL(1)^n} N)\\
    &\cong \scL(n + 1) \ltimes_{\scL(1)^{n + 1}} ((\scL(m) \times \scL(1)^n) \ltimes_{\scL(1)^{m + n}} (M \wedge N)).
\end{align*} But by \Cref{lem:hopkins} we have
\begin{align*}
    \scL(m + n) \ltimes_{\scL(1)^{m + n}} (M \wedge N) &\cong (\scL(n + 1) \times_{\scL(1)^{n + 1}} \scL(m) \times \scL(1)^n) \ltimes_{\scL(1)^{m + n}} (M \wedge N)\\
    &\cong \scL(n + 1) \ltimes_{\scL(1)^{n + 1}} ((\scL(m) \times \scL(1)^n)\ltimes_{\scL(1)^{m + n}} (M \wedge N)),
\end{align*} so we conclude that
\begin{equation*}
    \Phi_{n + 1}\mu(\Phi_m M, N) \cong \scL(m + n) \ltimes_{\scL(1)^{m + n}} (M \wedge N) = \Phi_{m + n}\mu(M, N),
\end{equation*} which is the desired isomorphism.
\end{proof}

A priori the internalized product $\wedge_\scL := \Phi_2\mu_{1, 1}$ of $\LL$-spectra is associative (see Theorem I.5.5 of \cite{EKMM}). One would hope that $S$ would be the unit of this product, but sadly this is not the case. The problem, as hinted at above, is the absence of absorbers in degree 0. As a consolation, we do, however, get a canonical map
\begin{equation*}
    S \wedge_\scL M \to M
\end{equation*} for any $\LL$-spectrum $M$, described in Proposition I.8.3 of \cite{EKMM}. In fact, this generalizes to higher degree spectra and to spaces besides $S^0$.

\begin{lem} \label{lem:Smodule_morphism}
Let $(M, \xi_M)$ be an $\LL^m$-spectrum, and let $X$ be a based topological space. Then there is a natural morphism of $\LL$-spectra
\begin{equation*}
    \Phi_{m+1}\mu(\Phi_0 X, M) \to \Phi_m\mu(X, M).
\end{equation*}
\end{lem}
\begin{proof}
    First, suppose $M$ is a free $\LL^m$-algebra $\LL^m E = \scL(1)^m \ltimes E$, where $E \in \scS[U^m]$. Then we seek a natural morphism
    \begin{equation*}
        \scL(m + 1) \ltimes_{\scL(1)^{m + 1}}(\Sigma^\infty X \wedge \scL(1)^m \ltimes E) \to \scL(m) \ltimes_{\scL(1)^m} (X \wedge \scL(1)^m \ltimes E).
    \end{equation*} First, we observe that $\gamma$ defines a morphism
    \begin{equation*}
        \scL(m + 1) \ltimes_{\scL(1)^{m + 1}} \scL(0) \times \scL(1)^m \to \scL(m),
    \end{equation*} since it coequalizes the pair of morphisms
    \begin{equation*}
        \begin{tikzcd}
            \scL(m + 1) \times \scL(1)^{m + 1} \times \scL(0) \times \scL(1)^m \ar[r, shift left = 0.35em, "\id \times \gamma^{m + 1}"] \ar[r, shift right = 0.35em, "\gamma \times \id"'] & \scL(m + 1) \times \scL(0) \times \scL(1)^m \period
        \end{tikzcd}
    \end{equation*} This yields a map from
    \begin{align*}
        \scL(m + 1) \ltimes_{\scL(1)^{m + 1}}(\Sigma^\infty X \wedge \scL(1)^m E) &= \scL(m + 1) \ltimes_{\scL(1)^{m + 1}}(\scL(0) \ltimes X \wedge \scL(1)^m \ltimes E)\\
        &\cong (\scL(m + 1) \times_{\scL(1)^{m + 1}} \scL(0) \times \scL(1)^m) \ltimes (X \wedge E)
    \end{align*} to
    \begin{equation*}
        \scL(m) \ltimes (X \wedge E).
    \end{equation*} But by \Cref{lem:hopkins} we have
    \begin{align*}
        \scL(m) \cong \scL(m) \times_{\scL(1)^m} \scL(1)^m,
    \end{align*} so
    \begin{align*}
        \scL(m) \ltimes (X \wedge E) &\cong \scL(m) \ltimes_{\scL(1)^m} (\scL(1)^m \ltimes (X \wedge E))\\
        &\cong \scL(m) \ltimes_{\scL(1)^m} (X \wedge (\scL(1)^m \ltimes E)).
    \end{align*} Composing the morphisms above yields the desired map.

    Now, if $M$ is any $\LL^m$-spectrum then there is a split coequalizer diagram
    \begin{equation*}
        \begin{tikzcd}
            \LL^m \LL^m M \ar[r, shift left = 0.25em] \ar[r, shift right = 0.25em] & \LL^m M \ar[r, two heads] & M \period
        \end{tikzcd}
    \end{equation*} Since split coequalizers are absolute colimits, the canonical morphisms corresponding to the free algebras in the diagram induce a morphism of coequalizer diagrams
    \begin{equation*}
        \begin{tikzcd}
            \Phi_{m + 1}\mu(\Phi_0 X, \LL^m \LL^m M) \ar[r, shift left = 0.25em] \ar[r, shift right = 0.25em] \ar[d] & \Phi_{m + 1}\mu(\Phi_0 X, \LL^m M) \ar[r, two heads] \ar[d] & \Phi_{m + 1}\mu(\Phi_0 X, M) \ar[d, dashed]\\
            \Phi_m\mu(X, \LL^m \LL^m M) \ar[r, shift left = 0.25em] \ar[r, shift right = 0.25em] & \Phi_m\mu(X, \LL^m M) \ar[r, two heads] & \Phi_m\mu(X, M)
        \end{tikzcd}
    \end{equation*}
    which induces the desired map for $M$.
\end{proof}

In general, the map in \Cref{lem:Smodule_morphism} need not be an isomorphism. To obtain the desired internalized category, we restrict our domain to only those $\LL_m$-spectra for which it is.

\begin{defn}\label{def:graded_Smodule}
    Given $m \geq 1$, let $\scM_m$, \emph{the category of $S$-modules of degree $m$}, be the full subcategory of $\scS[\LL]_m$ spanned by those $\LL^m$-algebras $M$ such that the canonical map
    \begin{equation*}
        \Phi_{m+1}\mu(\Phi_0 X, M) \to \Phi_m\mu(X, M)
    \end{equation*} of \Cref{lem:Smodule_morphism} an isomorphism.

    Additionally, let $\scM_0 := \Top_*$, the category of based spaces, and let $\Phi_0 := \scL(0) \ltimes (-) = \Sigma^\infty$ (which, by Proposition II.1.2 of \cite{EKMM}, defines a functor from based spaces to $S$-modules). Then, by way of the following proposition, we can naturally extend $\mu$ to a graded monoidal product on $\scM_\sbull$, with unit $S^0$.
\end{defn}

\begin{prop} \label{prop:smash_of_module}
    Let $M \in \scM_m$, and let $N \in \scM_n$. Then $\mu(M, N) \in \scM_{m + n}$.
\end{prop}
\begin{proof} 
    Let $X$ be a based topological space. We aim to show that the map
    \begin{equation*}
        \Phi_{m + n + 1}\mu(\Phi_0 X, \mu(M, N)) \to \Phi_{m + n}\mu(X, \mu(M, N))
    \end{equation*} is an isomorphism.
    
    First, suppose $m, n > 0$. From the associativity of $\mu$, we have
    \begin{equation*}
        \Phi_{m + n + 1}\mu(\Phi_0 X, \mu(M, N)) \cong \Phi_{m + n + 1}\mu(\mu(\Phi_0 X, M), N).
    \end{equation*} Now, since $\Phi_0 X$ and $\mu(M, N)$ both have positive degree, we can apply \Cref{lem:higher_absorbers} to conclude that
    \begin{equation*}
        \Phi_{m + n + 1}\mu(\mu(\Phi_0 X, M), N) \cong \Phi_{n + 1}\mu(\Phi_{m + 1}\mu(\Phi_0 X, M), N).
    \end{equation*} But $M$ is a degree $m$ $S$-module, so we have
    \begin{equation*}
        \Phi_{n + 1}\mu(\Phi_{m + 1}\mu(\Phi_0 X, M), N) \cong \Phi_{n + 1}\mu(\Phi_m\mu(X, M), N).
    \end{equation*} Since $N$ has positive degree, we can again utilize \Cref{lem:higher_absorbers}, from which we obtain
    \begin{equation*}
        \Phi_{n + 1}\mu(\Phi_m\mu(X, M), N) \cong \Phi_{m + n}\mu(\mu(X, M), N).
    \end{equation*} Finally, we again use the associativity of $\mu$ to conclude that
    \begin{equation*}
        \Phi_{m + n + 1}\mu(\Phi_0 X, \mu(M, N)) \cong \Phi_{m + n}\mu(X, \mu(M, N)),
    \end{equation*} as desired.

    Next, suppose (without loss of generality) that $m = 0$ (i.e. that $M$ is a based topological space) and $n > 0$. (The case when $m > 0$ and $n = 0$ follows from the symmetry of $\mu$.) Using the associativity of $\mu$ we have 
    \begin{equation*}
        \Phi_{n + 1}\mu(\Phi_0 X, \mu(M, N)) \cong \Phi_{n + 1}\mu(\mu(\Phi_0 X, M), N).
    \end{equation*} Moreover, we have
    \begin{equation*}
        \mu(\Phi_0 X, M) = (\Sigma^\infty X) \wedge M \cong (\scL(0) \ltimes X) \wedge M,
    \end{equation*} and by properties of the twisted half-smash product it follows that
    \begin{equation*}
        \mu(\Phi_0 X, M) \cong \scL(0) \ltimes (X \wedge M) \cong \Sigma^\infty(X \wedge M) = \Phi_0\mu(X, M).
    \end{equation*} Therefore, we have
    \begin{equation*}
        \Phi_{n + 1}\mu(\mu(\Phi_0 X, M), N) \cong \Phi_{n + 1}\mu(\Phi_0\mu(X, M), N),
    \end{equation*} and since $N$ is a degree $n$ $S$-module
    \begin{equation*}
        \Phi_{n + 1}\mu(\Phi_0\mu(X, M), N) \cong \Phi_n(\mu(X, M), N).
    \end{equation*} Applying the associativity of $\mu$ again yields the isomorphism
    \begin{equation*}
        \Phi_{n + 1}\mu(\Phi_0 X, \mu(M, N)) \cong \Phi_n(X, \mu(M, N)).
    \end{equation*}

    Lastly, if both $m$ and $n$ are $0$ then $M$ and $N$ are based topological spaces, in which case so is $\mu(M, N)$. Thus $\mu(M, N)$ is likewise an $S$-module of degree $0$, completing the proof.
\end{proof}

It follows that $\mu$ restricts to a graded monoidal product on $\scM_\sbull$.

\begin{rem}
    Note that $\scM_1 = \scM_S$, the category of $S$-modules (see Definition II.1.1 of \cite{EKMM}). Certainly each spectrum in $\scM_1$ is an $S$-module, by taking $X = S^0$. Conversely, every classical $S$-module is a degree $1$ $S$-module by Proposition II.1.4 of \cite{EKMM}.
\end{rem}

\begin{rem}
    If $M \in \scM_m$ then $\Phi_m(M) \in \scM_1$, so $\Phi_m$ restricts to a functor $\scM_m \to \scM_1$. This is because, given a based space $X$, we have isomorphisms
    \begin{align*}
        \Phi_{2}\mu(\Phi_0 X, \Phi_m M) \cong \Phi_{m + 1}\mu(\Phi_0 X, M) \cong \Phi_m\mu(X, M),
    \end{align*} where the first map is the isomorphism described in \Cref{lem:higher_absorbers}, and the second follows from the fact that $M$ a degree $m$ $S$-module.
\end{rem}

We can now verify that the functors $\Phi_\sbull$ indeed define an internalization of $\mu$. For this, we will need a preliminary result, which generalizes Lemma I.5.4 of \cite{EKMM}. 

\begin{thm}\label{thm:wedgeS_internal}
    The functors $\Phi$ define an internalization on the graded monoidal category $(\scM_\sbull, \mu, S^0)$, and the internalized product on $\scM_S$ is precisely $\wedge_S$.
\end{thm}
\begin{proof} We aim to construct absorbers for $\Phi$, as defined in \Cref{def:internalization}. We will content ourselves with identifying the left absorbers, as the right absorbers are built through an identical procedure, and omit the verification of the coherence axioms. 

Let $M \in \scM_m$, $N \in \scM_n$. We seek a natural isomorphism
\begin{equation*}
    \Phi_{n + 1}\mu(\Phi_m M, N) \cong \Phi_{m + n}\mu(M, N).
\end{equation*} First, if $m, n > 0$ then the left (and right) absorbers are given by restricting the isomorphisms of \Cref{lem:higher_absorbers} to $S$-modules.

Next, suppose $m = 0$. Then $M$ is a based topological space, and since $N$ is a degree $n$ $S$-module we have a natural isomorphism
\begin{equation*}
    \Phi_{n + 1}\mu(\Phi_0 M, N) \cong \Phi_n\mu(M, N)
\end{equation*} by \Cref{def:graded_Smodule}.

Finally, suppose $n = 0$. Then $N$ is a based topological space, and thus we have natural isomorphisms
\begin{align*}
    \Phi_1 \mu(\Phi_m M, N) &= \scL(1) \ltimes_{\scL(1)} ((\scL(m) \ltimes_{\scL(1)^m} M) \wedge N)\\
    &\cong (\scL(m) \ltimes_{\scL(1)^m} M) \wedge N\\
    &\cong \scL(m) \ltimes_{\scL(1)^m} (M \wedge N)\\
    &= \Phi_m\mu(M, N),
\end{align*} completing the proof.
\end{proof}

Having shown that the smash product $\wedge_S$ of $S$-modules is an internalization of $\wedge$, we can immediately conclude that it is symmetric, using the formal properties of $\Phi$. 

\begin{lem} \label{lem:M_sym}
    $(\scM_\sbull, \mu, S^0) $ is a symmetric graded monoidal category, and $\Phi$ is a symmetric internalization.
\end{lem}
\begin{proof}
    The symmetry of $(\scM_\sbull, \mu, S^0)$ comes from that of the external smash product $\wedge$. Naturality of the isomorphisms which define the left and right absorbers in the proof of \Cref{thm:wedgeS_internal} shows that $\Phi$ is symmetric.
\end{proof}

\begin{prop}
    The monoidal product $\wedge_S$ on the category of $S$-modules $\scM_S$ is symmetric.
\end{prop}
\begin{proof}
    This follows immediately from \Cref{lem:M_sym} and \Cref{cor:symmetric_internalization}.
\end{proof}

\section*{Acknowledgments}

I am immensely grateful to Peter May for his mentorship and patience. Without his support I would never have had the opportunity to work on this project, and it was a great privilege to study under him. I'd also like to thank Peter for the colossal amount of time and effort he puts into organizing the University of Chicago REU, and, of course, for the invitation to participate in it during the Summer 2022 program.

I also owe a great deal to both David Nadler and Peter Haine for their support during my time at UC Berkeley. Their guidance throughout my undergraduate career meant the world to me, and I would not be where I am today without their help. It was a privilege to have had the chance to learn from both of them.

Lastly, I would like to thank Maximilien Péroux for very generously taking the time to help me finish this paper, and Nelson Niu for his advice on handling various aspects of monoidal coherence. I am also grateful to both of them for their input on further applications of my framework.

\nocite{*}

\newpage
\printbibliography

\end{document}